\documentclass[11pt]{amsart}

\usepackage{enumitem}
\usepackage{psfrag}
\usepackage{graphicx}
\usepackage{pinlabel}
\usepackage{graphicx}
\usepackage{amsmath}
\usepackage{amssymb}
\usepackage{amscd}
\usepackage{picinpar}
\usepackage{times}
\usepackage{pb-diagram}
\usepackage{graphicx}
\usepackage{wrapfig}
\usepackage{xspace}
\usepackage{hyperref}
\usepackage{url}
\usepackage{caption}
\usepackage{subcaption}
\usepackage{fancyhdr}
\newtheorem{theorem}{Theorem}[section]
\newtheorem{lemma}[theorem]{Lemma}
\newtheorem{proposition}[theorem]{Proposition}
\newtheorem{corollary}[theorem]{Corollary}
\newtheorem{definition}[theorem]{Definition}

\newtheorem{remark}[theorem]{Remark}
\newtheorem*{theorem*}{Theorem}
\def\C{\mathbb{C}}

\def\T{\mathcal{T}}

\usepackage{overpic}
\newcommand{\bal}{\begin{aligned}}      \newcommand{\eal}{\end{aligned}}
\newcommand{\ba}{\begin{array}}      \newcommand{\ea}{\end{array}}
\newcommand{\bc}{\begin{center}}     \newcommand{\ec}{\end{center}}
\newcommand{\be}{\begin{enumerate}}  \newcommand{\ee}{\end{enumerate}}
\newcommand{\beq}{\begin{eqnarray}}  \newcommand{\eeq}{\end{eqnarray}}
\newcommand{\beQ}{\begin{eqnarray*}} \newcommand{\eeQ}{\end{eqnarray*}}
\newcommand{\bi}{\begin{itemize}}    \newcommand{\ei}{\end{itemize}}
\newcommand{\bt}{\begin{tabular}}    \newcommand{\et}{\end{tabular}}
\newcommand{\bdm}{\begin{displaymath}} \newcommand{\edm}{\end{displaymath}}

\def\qed{\hfill{Q.E.D.}\smallskip}

\begin{document}

\title{\bf Rigidity of infinite inversive distance circle packings in the plane}
\author{Yanwen Luo, Xu Xu, Siqi Zhang}

\address{Department of Mathematics, Oklahoma State University, Stillwater
 OK, 74078}
 \email{yanwen.luo@okstate.edu}

\address{School of Mathematics and Statistics, Wuhan University, Wuhan, 430072, P.R.China}
 \email{xuxu2@whu.edu.cn}

\address{School of Mathematics and Statistics, Wuhan University, Wuhan, 430072, P.R.China}
 \email{2014301000156@whu.edu.cn}

\maketitle

\begin{abstract}
In 2004, Bowers-Stephenson \cite{BS} introduced the inversive distance circle packings as a natural generalization of
Thurston's circle packings.
They further conjectured the rigidity of infinite inversive distance circle packings in the plane.
Motivated by the recent work of Luo-Sun-Wu \cite{LSW} on Luo's vertex scaling,
we prove Bowers-Stephenson's conjecture for inversive distance circle packings in the hexagonal triangulated plane.
This generalizes Rodin-Sullivan's famous result \cite{RS} on the rigidity of infinite tangential circle packings in the hexagonal  triangulated plane.
The key tools include a maximal principle for generic weighted Delaunay inversive distance circle packings
and a ring lemma for the inversive distance circle packings in the hexagonal triangulated plane.
\end{abstract}



%\tableofcontents

\section{Introduction}
In 1985, Thurston \cite{Th2} conjectured that the Riemann mapping
for simply connected domains in the plane could be approximated by tangential circle packings.
Thurston's conjecture was solved elegantly by Rodin-Sullivan \cite{RS} by proving the rigidity
of infinite tangential circle packings in the hexagonal triangulated plane.
Since then, there have been lots of important works on the rigidity of infinite circle packings in the plane. See \cite{H2, Sch, St} and others.

Motivited by Thurston's circle packings \cite{Th}, Bowers-Stephenson \cite{BS} introduced
the inversive distance circle packings as a natural generalization.
They further conjectured the rigidity of infinite inversive distance circle packings in the plane.
In this paper, we prove Bowers-Stephenson's conjecture for weighted Delaunay inversive distance circle packings in
the hexagonal triangulated plane.
The proof is accomplished by establishing a maximal principle for generic weighted Delaunay inversive distance circle packings
and a ring lemma for inversive distance circle packings in the hexagonal triangulated plane.
The main idea comes from the recent work of Luo-Sun-Wu \cite{LSW}, in which the infinite rigidity of Luo's vertex scaling \cite{Luo CCM} in the
hexagonal plane was proved.

Suppose $S$ is a topological surface  and $\mathcal{T}$
is a triangulation of $S$.
We use $V = V(\mathcal{T})$, $E = E(\mathcal{T})$ and $F = F(\mathcal{T})$ to
denote the set of vertices, edges, and faces of $\mathcal{T}$ respectively.
A piecewise linear metric $d$ (PL metric for simplicity) on $(S, \mathcal{T})$
is a flat cone metric on $S$ such that each face in $F$ in the metric $d$ is isometric to a Euclidean triangle.
For simplicity,  a PL metric on $(S, \mathcal{T})$ is
represented as a function $l: E\rightarrow (0,+\infty)$ satisfying the strict triangle inequality.
For a PL metric $l: E\rightarrow (0,+\infty)$ on $(S, \mathcal{T})$, the combinatorial curvature is a map $K: V\rightarrow (-\infty, 2\pi)$
sending an interior vertex $v\in V$ to $2\pi$ minus the sum of angles at $v$ and
a boundary vertex $v\in V$ to $\pi$ minus the angles at $v$.
The combinatorial curvature $K$ on a compact triangulated surface $(S, \mathcal{T})$ satisfies the discrete Gauss-Bonnet formula \cite{CL}
\begin{equation}\label{discrete Gauss-Bonnet formula}
  \sum_{i\in V}K_i=2\pi\chi(S).
\end{equation}
A PL metric is called flat if $K(v)=0$ for any interior vertex $v$.

\begin{definition}[\cite{BS}]\label{discrere conformal for idcp}
Suppose $(S, \mathcal{T})$ is a triangulated surface with a weight $I: E\rightarrow (-1, +\infty)$.
  A PL metric $l: E\rightarrow (0,+\infty)$ on the weighted triangulated surface $(S, \mathcal{T}, I)$
  is an inversive distance circle packing metric
  if there exists a function $u: V \to \mathbb{R}$  such that for any edge $e\in E$ with end points $v$
  and $v'$, the length $l(e)$ is given by
\begin{equation}
\label{length1 introduction}
  l(e)=\sqrt{e^{2u(v)}+e^{2u(v')}+2 I(e)e^{u(v)+u(v')}}.
\end{equation}
The function $u: V \to \mathbb{R}$ is called a label on $(S, \mathcal{T}, I)$.
Two inversive distance circle packing metrics $l$ and $\tilde l$ on $(S, \mathcal{T}, I)$ are conformally equivalent.
In this case, we set $w=\tilde{u}-u$ and denote $\tilde{l}$ by $w*l$.  $w$ is called a discrete conformal factor.
\end{definition}

Thurston's circle packing metric \cite{Th} is a special type of
inversive distance circle packing metric with $I\in [0,1]$ in (\ref{length1 introduction}).
If we set $r(v)=e^{u(v)}$ and $r(v')=e^{u(v')}$, then the weight $I(e)$ in (\ref{length1 introduction}) is the inversive distance of
the two circles centered at $v$ and $v'$ with radii $r(v)$ and $r(v')$ respectively in the plane.
The map $r: V\rightarrow (0, +\infty)$ is also said to be an \textit{inversive distance circle packing} on the weighted
triangulated surface $(S, \mathcal{T}, I)$.
The rigidity of finite inversive distance circle packings on a weighted
triangulated closed surfaces $(S, \mathcal{T}, I)$ has been proved in \cite{Guo,Luo GT,Xu AIM,Xu MRL}.
The main focus of this paper is to provide an affirmative answer to Bowers-Stephenson conjecture on the rigidity
of infinite inversive distance circle packings in the hexagonal triangulated plane.
To state the main result, we need to introduce the following notions for inversive distance circle packings.

Assume that $r: V\rightarrow (0, +\infty)$ is an inversive distance circle packing on a weighted triangulated surface $(S, \mathcal{T}, I)$
with $I: E\rightarrow (-1, +\infty)$.
Let $\triangle v_1v_2v_3$ be a Euclidean triangle in the plane isometric to a face in $(S, \mathcal{T}, I, r)$,
each vertex $v_i$ of which is attached with a circle of radius $r_i=r(v_i)$ centered at the vertex.
The \textit{power distance} of a point $p$ in the plane to the vertex $v_i$ is defined to be $\pi_i(p)=d^2(v_i,p)-r_i^2$, where $d(v_i,p)$ is the Euclidean distance between $p$ and the vertex $v_i$. The \textit{geometric center} $C_{123}$ of $\triangle v_1v_2v_3$
is the unique point in the plane having the same power distance to the vertices $v_1, v_2, v_3$.
Denote $h_{jk,i}$ as the signed distance of the geometric center $C_{123}$ to the edge ${v_jv_k}$, which is positive if $C_{123}$ is in the same side of the line $v_jv_k$ as $\triangle v_1v_2v_3$ and negative otherwise.
Please refer to \cite{Glickenstein JDG, Glickenstein preprint, GT} for more information on the geometric centers of discrete conformal structures on manifolds.
\begin{definition}[\cite{Glickenstein DCG}]\label{defn of weighted Delaunay}
Suppose $r: V\rightarrow (0, +\infty)$ is an inversive distance circle packing on a weighted triangulated surface $(S, \mathcal{T}, I)$
with $I: E\rightarrow (-1, +\infty)$.
$v_iv_j\in E$ is an edge shared by two adjacent triangles $\triangle v_iv_jv_k$ and $\triangle v_iv_jv_m$ in $\mathcal{T}$.
The edge $v_iv_j$ is weighted Delaunay in $(S, \mathcal{T}, I, r)$ if
$$h_{ij,k}+h_{ij,m}\geq 0.$$
The triangulation $\mathcal{T}$ is weighted Delaunay in $r$ if every interior edge is weighed Delaunay.
\end{definition}
For simplicity, we call $r$ as a weighted Delaunay inversive distance circle packing on $(S, \mathcal{T}, I)$,
if the triangulation $\mathcal{T}$ is weighted Delaunay in $r$.
In this case, we also say that the PL metric induced by $r$ on $(S, \mathcal{T}, I)$ is weighted Delaunay.
There are other equivalent definitions for weighted Delaunay triangulations. Please refer to \cite{AK, Glickenstein DCG, Glickenstein preprint} and others.

The weight $I: E\to (-1, +\infty)$ on a triangulated surface $(S, \mathcal{T})$ is \emph{regular} if there is no adjacent triangles $t_1$ (with edges $a,b,e$) and $t_2$ (with edges $c,d,e$) in $F$ such that $I(e)=1, I(a)=-I(b), I(c)=-I(d).$
For a hexagonal triangulation $\mathcal{T}$ of the plane, we can take $V$ as the lattice $L=\{m\vec{v}_1+n\vec{v}_2| m,n\in \mathbb{Z}, \vec{v}_1=1, \vec{v}_2=e^{i\frac{\pi}{3}}\}$, in which the addition of vertices could be defined.
A weight $I: E\to (-1, +\infty)$ on the hexagonal triangulated plane is \textit{translating invariant} if $I(e+\delta)=I(e)$ for any $e\in E, \delta\in L$,
where $e+\delta$ is an edge with end points $v+\delta$ and $v'+\delta$ if the edge $e\in E$ has end points $v$ and $v'$.

The main result of this paper is as follows.

\begin{theorem}
\label{infrigidity introduction}
Let $(\mathbb{C}, \mathcal{T}_{st})$ be a hexagonal triangulated plane.
$I$ is a regular, translating invariant weight defined on the edges with values in $(-\frac{1}{2}, 1]$ or $[0,+\infty)$ and satisfying the following structure condition
		\begin{equation}
		\label{structure condition}
			I(e_i)+I(e_j)I(e_k)\geq 0, \{i,j,k\}=\{1,2,3\}
		\end{equation}
for any triangle in $\mathcal{T}$ with edges $e_1, e_2, e_3$.
Assume $l$ is a weighted Delaunay inversive distance circle packing metric on $(\mathbb{C}, \mathcal{T}_{st}, I)$
induced by a constant label.
If $(\mathbb{C}, \mathcal{T}_{st}, I, w*l)$ is a weighted Delaunay triangulated surface isometric to an open set in the plane, then $w$ is a constant function.
\end{theorem}

Theorem \ref{infrigidity introduction} generalizes Rodin-Sullivan's famous result \cite{RS} on the rigidity of infinite tangential circle packings in the hexagonal plane, which corresponds to $I\equiv 1$.

The paper is organized as follows.
In Section \ref{section 2}, we give some preliminaries on the inversive distance circle packings
and weighted Delaunay triangulations.
In Section \ref{section 3}, we derive the maximal principle for generic inversive distance circle packings
 and the ring lemma for inversive distance circle packings in the hexagonal triangulated plane.
We also study the properties of inversive distance circle packings on spiral hexagonal triangulations in this section.
In Section \ref{section 4}, we prove a generalized version of Theorem \ref{infrigidity introduction}, i.e. the rigidity of infinite inversive distance circle packings in the hexagonal triangulated plane.
\\
\\
\textbf{Acknowledgements}\\[8pt]
The research of the second author is supported by National Natural Science Foundation of China
under grant no. 12471057.

\section{Preliminaries on inversive distance circle packings and weighted Delaunay triangulations}\label{section 2}
%In this section, we collect some basic properties of inversive distance circle packings and weighted Delaunay triangulations.

Let $(S, \mathcal{T}, I)$ be a weighted triangulated surface with the weight $I: E\rightarrow (-1, +\infty)$.
We use $v_i$ to denote a vertex in $V$, $v_iv_j$ to denote an edge in $E$ and $\triangle v_iv_jv_k$ to denote a face in $F$.
We further denote $f_i = f(v_i)$ if $f$ is a function defined on $V$,  $f_{ij} = f(v_iv_j)$ if $f$ is a function defined on $E$, and $f_{ijk} = f(\triangle v_iv_jv_k)$ if $f$ is a function defined on $F$.

\subsection{Basic properties of inversive distance circle packings}
For any function $u: V\rightarrow \mathbb{R}$, the formula (\ref{length1 introduction})
gives a positive number $l(e)$ for any edge $e\in E$ since $I(e)>-1$.
However, for a face $\triangle v_iv_jv_k$ in $(S, \mathcal{T}, I)$, the positive numbers $l_{ij}, l_{ik}, l_{jk}$  may not satisfy
the \textit{strict triangle inequality}
\begin{equation}\label{strict triangle inequality}
l_{rs}< l_{rt}+l_{st}, \{r,s,t\}=\{i,j,k\}.
\end{equation}
The label $u: V\rightarrow \mathbb{R}$ is said to be \textit{admissible}
if the function $l: E\rightarrow (0, +\infty)$ determined by $u: V\rightarrow \mathbb{R}$
via the formula (\ref{length1 introduction}) satisfies
the strict triangle inequality (\ref{strict triangle inequality}) for every face in $(S, \mathcal{T}, I)$, i.e. $l: E\rightarrow (0, +\infty)$ is
a PL metric on $(S, \mathcal{T})$.
We also say that the corresponding inversive distance circle packing $r:V\rightarrow (0, +\infty)$
on $(S, \mathcal{T}, I)$ with $r_i=e^{u_i}$ is admissible, if it causes no confusion in the context.
The admissible space of inversive distance circle packings on $(S, \mathcal{T}, I)$ is the set of admissible inversive distance circle packings on
$(S, \mathcal{T}, I)$.
For an admissible inversive distance circle packing $r$ on $(S, \mathcal{T}, I)$, every face in $(S, \mathcal{T}, I)$ is isometric to a \textit{non-degenerate (Euclidean) triangle}
with edge lengths given by (\ref{length1 introduction}).
We also say that $r: V\rightarrow (0,+\infty)$ generates a PL metric on $(S, \mathcal{T}, I)$ for simplicity in this case.

If three positive numbers $l_{ij}, l_{ik}, l_{jk}$ satisfy the \textit{triangle inequality }
\begin{equation}\label{triangle inequality}
l_{rs}\leq l_{rt}+l_{st}, \{r,s,t\}=\{i,j,k\},
\end{equation}
then $l_{ij}, l_{ik}, l_{jk}$ generates a \textit{generalized (Euclidean) triangle} $\triangle v_iv_jv_k$.
If $l_{ij}=l_{ik}+l_{jk}$, the generalized triangle $\triangle v_iv_jv_k$ is flat at $v_k$,
the inner angle at which is defined to be $\pi$. In this case, the generalized triangle is referred as a \textit{degenerate triangle}.
A function $l: E\rightarrow (0,+\infty)$ is called a \textit{generalized PL metric} on $(S,\mathcal{T})$
if the triangle inequality (\ref{triangle inequality}) is satisfied for
every face in $(S,\mathcal{T})$.
The PL metric is a special type of generalized PL metric.
The combinatorial curvature for generalized PL metrics is defined the same as the PL metrics
and still satisfies the discrete Gauss-Bonnet formula (\ref{discrete Gauss-Bonnet formula}) on compact triangulated surfaces.
A generalized PL metric $l:E\rightarrow (0,+\infty)$ is called a \textit{generalized inversive distance circle packing metric} on a weighted
triangulated surface $(S, \mathcal{T}, I)$ if there exists a map $u: V\rightarrow \mathbb{R}$ such that $l$ is determined by $u$
via the formula (\ref{length1 introduction}).
In this case, the map $r: V\rightarrow (0,+\infty)$ with $r_i=e^{u_i}$ is said to be
\textit{a generalized inversive distance circle packing} on $(S, \mathcal{T}, I)$.
	
	\begin{lemma}[\cite{Guo, Xu AIM, Xu MRL}]\label{basic property I of IDCP}
	Let $\triangle v_1v_2v_3$ be a face in $(S, \mathcal{T})$ with three weights $I_1,I_2,I_3\in (-1, +\infty)$ defined on edges opposite to the vertices $v_1, v_2, v_3$ respectively.
$u: \{v_1, v_2, v_3\}\rightarrow \mathbb{R}$ is a function defined on the vertices and the edge lengthes
are defined by
		\begin{equation}
		\label{length2}
			l_{ij}= \sqrt{e^{2u_i}+e^{2u_j}+2e^{u_i+u_j}I_{k}} =  \sqrt{r_i^2+r_j^2+2r_ir_jI_{k}},
		\end{equation}
where $r_i=e^{u_i}$, $\{i,j,k\}=\{1,2,3\}$.
\begin{description}
  \item[(a)] $l_{12}, l_{13}, l_{23}$ generate a non-degenerate Euclidean triangle if and only if $Q>0$, where
	\begin{equation}\label{defn of Q}
\begin{aligned}
		Q=\kappa_1^2(1-I^2_{1})+\kappa_2^2(1-I^2_{2})+\kappa_3^2(1-I^2_{3})
		+2\kappa_1\kappa_2\gamma_{3}+2\kappa_1\kappa_3\gamma_{2}+2\kappa_2\kappa_3\gamma_{1}
\end{aligned}
	\end{equation}
with $\gamma_{i}:=I_{i}+I_{j}I_{k}, \ \kappa_i:=r_i^{-1}$.
As a result,  $l_{12}, l_{13}, l_{23}$ generate a degenerate Euclidean triangle if and only if $Q=0$.
Specially, if the weights $I_1, I_2, I_3\in (-1, 1]$ satisfy the structure condition (\ref{structure condition}),  i.e. $\gamma_i=I_i+I_jI_k\geq 0, \{i,j,k\}=\{1,2,3\}$, then $l_{12}, l_{13}, l_{23}$ generate a non-degenerate Euclidean triangle for any $(u_1,u_2,u_3)\in \mathbb{R}^3$.
  \item[(b)] Assume that the weights $I_1, I_2, I_3\in (-1, +\infty)$ satisfy the structure condition (\ref{structure condition}),
and $u: \{v_1, v_2, v_3\}$
$\rightarrow \mathbb{R}$ generates a non-degenerate
Euclidean triangle $\triangle v_1v_2v_3$ with the edge length given by the formula (\ref{length2}).
Let $\theta_i$ be its inner angle at $v_i$. Then
	    $$\frac{\partial \theta_i}{\partial u_j}=\frac{\partial \theta_j}{\partial u_i}=\frac{h_{ij,k}}{l_{ij}}, \ \ \
\frac{\partial \theta_i}{\partial u_i}=-\frac{\partial \theta_i}{\partial u_j}-\frac{\partial \theta_i}{\partial u_k}<0, $$
		where
		\begin{equation}\label{h_ij,k}
			\begin{aligned}
				h_{ij,k}
			=\frac{r_1^2r_2^2r_3^2}{A_{123}l_{ij}}[\kappa_k^2(1-I_k^2)+\kappa_j\kappa_k\gamma_{i}+\kappa_i\kappa_k\gamma_{j}]
				=\frac{r_1^2r_2^2r_3^2}{A_{123}l_{ij}}\kappa_kh_k
			\end{aligned}
		\end{equation}
		with
$A_{123}=l_{12}l_{13}\sin\theta_1$ and
		\begin{equation}
		\label{h_i}
			\begin{aligned}
h_k=\kappa_k(1-I_k^2)+\kappa_i\gamma_{j}+\kappa_j\gamma_{i}.
			\end{aligned}
		\end{equation}
	Moreover, under the structure condition (\ref{structure condition}),
the Jacobian matrix $\Lambda_{123}=\frac{\partial(\theta_1, \theta_2, \theta_3)}{\partial (u_1, u_2, u_3)}$ for admissible $u$
	is negative semi-definite with one dimensional kernel $\{t(1,1,1)|t\in \mathbb{R}\}$.
  \item[(c)] Under the structure condition (\ref{structure condition}), if $(u_1,u_2,u_3)\in \mathbb{R}^3$ is not admissible,
then one of $h_1, h_2, h_3$ is negative and the other two are positive.
Specially, if $(u_1,u_2,u_3)\in \mathbb{R}^3$ generates a degenerate triangle
$\triangle v_1v_2v_3$ having $v_3$ as the flat vertex, then $h_1>0, h_2>0, h_3<0$ at $(u_1,u_2,u_3)$, which further implies $I_3>1$
by $h_3<0$.
  \item[(d)] If the structure condition (\ref{structure condition}) is satisfied and there exists $i\in \{1,2,3\}$ such that $I_i>1$, then      $$\Delta_{123}:=I_{1}^2+I_{2}^2+I_{3}^2+2I_1I_2I_3-1>0.$$
      As a result, if $\Delta_{123}\leq 0$ and the structure condition (\ref{structure condition}) is satisfied, then $I_1, I_2, I_3\in (-1, 1]$, $Q(r)>0$ for any $r\in \mathbb{R}^3_{>0}$, and the triangle $\triangle v_1v_2v_3$ generated by any $r\in \mathbb{R}^3_{>0}$ is always non-degenerate.
\end{description}
	\end{lemma}

%The admissible space of labels on a face in $(S, \mathcal{T}, I)$ is characterized in the following result.
As an application of Lemma \ref{basic property I of IDCP}, we have the following characterization of
the admissible space of inversive distance circle packings on a weighted triangle and extension of inner angles
for generalized triangles generated by inversive distance circle packings.
	
\begin{lemma}[\cite{Xu MRL}]\label{simply connect of admi space with weight}%[\cite{Guo,X2}]
Suppose $\triangle v_1v_2v_3$ is a face in $(S, \mathcal{T}, I)$ with the weight $I: E\rightarrow (-1, +\infty)$
satisfying the structure condition (\ref{structure condition}).  Then the admissible
space $\Omega_{123}$ of inversive distance circle packings $(r_1, r_2, r_3)\in \mathbb{R}^3_{>0}$
on $\triangle v_1v_2v_3$ is
$$\Omega_{123}=\mathbb{R}^3_{>0}\setminus \sqcup_{i\in P}V_i,$$
		where $P=\{i\in \{1,2,3\}|I_i>1\}$, $\sqcup_{i\in P}V_i$ is a disjoint union of
		$$V_i=\left\{(r_1, r_2, r_3)\in \mathbb{R}^3_{>0}|\kappa_i\geq \frac{-B_i+\sqrt{\Delta_i}}{2A_i}\right\}$$
		with
		\begin{equation}
			\begin{aligned}
				A_i=&I^2_{i}-1,\\
				B_i=&-2(\kappa_j\gamma_{k}+\kappa_{k}\gamma_j)\leq0,\\
				\Delta_i
				=&4(I_1^2+I_2^2+I_{3}^2+2I_1I_2I_{3}-1)(\kappa_j^2+\kappa_{k}^2+2\kappa_j\kappa_{k}I_i).
			\end{aligned}
		\end{equation}
As a result, $\Omega_{123}$ is nonempty and simply connected with analytic boundary.
Furthermore, the inner angles of $\triangle v_1v_2v_3$ could be uniquely continuously extended by constants as follows
	\begin{equation*}
		\begin{aligned}
			\widetilde{\theta}_i(r_1,r_2,r_3)=\left\{
			\begin{array}{ll}
				\theta_i, & \hbox{if $(r_1,r_2,r_3)\in \Omega_{123}$;} \\
				\pi, & \hbox{if $(r_1,r_2,r_3)\in V_i$;} \\
				0, & \hbox{otherwise.}
			\end{array}
			\right.
		\end{aligned}
	\end{equation*}
As a corollary, if $v_i$ is the flat point of the degenerate triangle $\triangle v_1v_2v_3$ generated by $(r_1,r_2,r_3)\in \mathbb{R}^3_{>0}$, then
$(r_1,r_2,r_3)\in \partial V_i$, i.e.
$$\kappa_i=\frac{-B_i+\sqrt{\Delta_i}}{2A_i}.$$
	\end{lemma}

\proof
We just need to prove the last part of the lemma, the other parts of the lemma have been proved in \cite{Xu MRL}.
By the first part of this lemma,
$(r_1,r_2,r_3)\in \partial \Omega_{123}$ in $\mathbb{R}^3_{>0}$, which is the disjoint union of $\partial V_1$, $\partial V_2$ and $\partial V_3$ in $\mathbb{R}^3_{>0}$.
Note that the inner angle of the degenerate triangle $\triangle v_1v_2v_3$ at $v_i$ is $\pi$ by assumption.
By the unique continuous extension of inner angles in the second part of this lemma, we have $(r_1,r_2,r_3)\in \partial V_i$, which implies
$\kappa_i=\frac{-B_i+\sqrt{\Delta_i}}{2A_i}.$
\qed

We prove the following results on inversive distance circle packings following Luo-Sun-Wu \cite{LSW}.
	
\begin{lemma}
		\label{interval}
Let $\triangle v_1v_2v_3$ be a face in $(S, \mathcal{T}, I)$ with the weight $I: E\rightarrow (-1, +\infty)$
satisfying the structure condition (\ref{structure condition}).
\begin{description}
  \item[(a)] For any fixed $r_i, r_j\in (0, +\infty)$, the set of $r_k\in (0, +\infty)$ such that $(r_i, r_j, r_k)$ is an admissible inversive distance
		circle packing on $\triangle v_1v_2v_3$  is an open interval.
As a result, if $ (r_i, r_j, \hat{r}_k)$ and $ (r_i, r_j, \bar{r}_k)$ are two generalized inversive distance circle packings on $\triangle v_1v_2v_3$ with $\hat{r}_k<\bar{r}_k$, then for any $r_k\in (\hat{r}_k,\bar{r}_k)$,
$ (r_i, r_j, r_k)$ generates a non-degenerate triangle $\triangle v_1v_2v_3$.
  \item[(b)] If $\triangle v_1v_2v_3$ generated by $(r_1,r_2,r_3)\in \mathbb{R}^3_{>0}$ is a degenerate triangle having $v_3$ as the flat vertex,  %$t>0$ the triangle generated by
% 		$(r_1, r_2, r_3+t)$ is a non-degenerate triangle. Furthermore, the quantity $h_{12,3}(t)$ is strictly increasing in $t$
	    then there exists $\epsilon>0$ such that $(r_1,r_2,r_3+t)\in \Omega_{123}$  and
		$\frac{\partial h_{12,3}}{\partial r_3}(r_1, r_2, r_3+t)>0$
for $t\in (0, \epsilon)$.
%{\color{red}Furthermore, there exists $\epsilon>0$ such that
%$$\frac{\partial h_{13,2}}{\partial r_3}(r_1, r_2, r_3+t)<0, \frac{\partial h_{23,1}}{\partial r_3}(r_1, r_2, r_3+t)<0$$
%for $t\in (0, \epsilon)$.}
\end{description}
	\end{lemma}
	
	\begin{proof}
    To prove part (a),
	without loss of generality, set $\{i,j\}=\{2,3\}$, $k=1$ and
$$
f(\kappa_1)=(1-I^2_{1})\kappa_1^2+2\kappa_1(\kappa_2\gamma_{3}+\kappa_3\gamma_{2})+
	\kappa_2^2(1-I^2_{2})+\kappa_3^2(1-I^2_{3})+2\kappa_2\kappa_3\gamma_{1}.
$$
By Lemma \ref{basic property I of IDCP} (a), we just need to show that the solution of $f(\kappa_1)>0$ with $\kappa_1\in (0, +\infty)$ is an open interval, which is included in the following three cases.
	
Case 1: If $I_1=1$, $f(\kappa_1)>0$ is equivalent to
	\begin{equation}\label{f(kappa1)>0 with I_1=1}
f(\kappa_1)=2\kappa_1(\kappa_2\gamma_{3}+\kappa_3\gamma_{2})+
	\kappa_2^2(1-I^2_{2})+\kappa_3^2(1-I^2_{3})+2\kappa_2\kappa_3\gamma_{1}>0.
\end{equation}
	If $\gamma_2=\gamma_3=0$, then $\gamma_2+\gamma_3=(1+I_1)(I_2+I_3)=0$, which implies $I_2+I_3=0$ by $I_1>-1$.
Therefore, $I_2, I_3\in (-1, 1)$ by $I_2, I_3\in (-1, +\infty)$, which further implies
	$f(\kappa_1)=
	\kappa_2^2(1-I^2_{2})+\kappa_3^2(1-I^2_{3})+2\kappa_2\kappa_3\gamma_{1}>0$
	for any $\kappa_2, \kappa_3\in (0, +\infty)$ in this case. Therefore, the solution of $f(\kappa_1)>0$ is $\mathbb{R}_{>0}$ in this case. \\
	If one of $\gamma_2$ and $\gamma_3$ is positive, then $\kappa_2\gamma_{3}+\kappa_3\gamma_{2}>0$ by the structure condition (\ref{structure condition}), which implies the solution of (\ref{f(kappa1)>0 with I_1=1}) is
	\begin{equation*} \kappa_1>-\frac{\kappa_2^2(1-I^2_{2})+\kappa_3^2(1-I^2_{3})+2\kappa_2\kappa_3\gamma_{1}}{2(\kappa_2\gamma_{3}+\kappa_3\gamma_{2})}.
	\end{equation*}
	This implies that the solution of $f(\kappa_1)>0$ with $\kappa_1>0$ is an open interval in this case.
	
Case 2: If $I_1\in (-1,1)$, then $1-I_1^2>0$ and
% 	\begin{equation}
% 		\begin{aligned}
% 			\Delta_1
% 			=4(I_1^2+I_2^2+I_3^2+2I_1I_2I_3-1)(\kappa_2^2+\kappa_3^2+2\kappa_2\kappa_3I_1).
% 		\end{aligned}
% 	\end{equation}
% 	The sign of $\Delta_1$ is determined by $I_1^2+I_2^2+I_3^2+2I_1I_2I_3-1$.
% 	If $I_1^2+I_2^2+I_3^2+2I_1I_2I_3-1<0$, then for any $\kappa_1\in (0, +\infty)$, we have $f(\kappa_1)>0$. If $I_1^2+I_2^2+I_3^2+2I_1I_2I_3-1\geq 0$, then by
	\begin{equation*}
		-\frac{b}{2a}=-\frac{\kappa_2\gamma_{3}+\kappa_3\gamma_{2}}{1-I_1^2}\leq 0,
	\end{equation*}
	which implies that the solution of the quadratic inequality $f(\kappa_1)>0$ with $\kappa_1>0$ is an open interval in $(0,+\infty)$ in this case.
	
Case 3: If $I_1\in (1, +\infty)$, then $f(\kappa_1)>0$ is equivalent to the following quadratic inequality
	$$(I^2_{1}-1)\kappa_1^2-2\kappa_1(\kappa_2\gamma_{3}+\kappa_3\gamma_{2})-
	\kappa_2^2(1-I^2_{2})-\kappa_3^2(1-I^2_{3})-2\kappa_2\kappa_3\gamma_{1}<0.$$
    In this case,
	\begin{equation*}
		-\frac{b}{2a}=\frac{\kappa_2\gamma_{3}+\kappa_3\gamma_{2}}{I_1^2-1}\geq 0,
	\end{equation*}
	and  the discriminant
	\begin{equation*}
		\begin{aligned}
			\Delta
			=4(I_1^2+I_2^2+I_3^2+2I_1I_2I_3-1)(\kappa_2^2+\kappa_3^2+2\kappa_2\kappa_3I_1)>0
		\end{aligned}
	\end{equation*}
by Lemma \ref{basic property I of IDCP} (d).
    This implies that the solution of $f(\kappa_1)>0$ with $\kappa_1>0$ is an open interval in $(0, +\infty)$ in this case.	

% 	Part (2) is based on the result in ??? \textcolor{red}{here need reference.}we have
% 	$$\frac{\partial \theta_1}{\partial u_1}=(1,0,0)\frac{\partial(\theta_1, \theta_2, \theta_3)}{\partial (u_1, u_2, u_3)}(1,0,0)^T<0$$
% 	by the fact that $\Lambda^E_{123}(\eta)=\frac{\partial(\theta_1, \theta_2, \theta_3)}{\partial (u_1, u_2, u_3)}$
% 	is negative semi-definite with one dimensional kernel $\{t(1,1,1)|t\in \mathbb{R}\}$.
% 	The other relationships in the Lemma has been proved in \cite{XX1}.\\
	
	To prove part (b), recall that the triangle $\triangle v_1v_2v_3$ is degenerate if and only if $Q=0$ by Lemma \ref{basic property I of IDCP} (a), where $Q$ is defined by (\ref{defn of Q}).
	By direct calculations, we have
	$\frac{\partial Q}{\partial \kappa_3}=2h_3<0$ at $(r_1,r_2,r_3)$
	by Lemma \ref{basic property I of IDCP} (c),
	which implies that
	$\frac{\partial Q}{\partial r_3}=\frac{\partial Q}{\partial \kappa_3}\frac{\partial \kappa_3}{\partial r_3}=-\frac{1}{r_3^2}\frac{\partial Q}{\partial \kappa_3}>0$
	around $(r_1, r_2, r_3)$.
	Therefore, for small $t>0$, $Q(r_1,r_2,r_3+t)>0$ and $(r_1,r_2,r_3+t)$ generates a non-degenerate triangle.
This can also be taken as a corollary of Lemma \ref{basic property I of IDCP} (c) and
Lemma \ref{simply connect of admi space with weight}.
	
	Recall that for a non-degenerate inversive distance circle packing on $\triangle v_1v_2v_3$, we have
$
h_{12,3}=\frac{r_1^2r_2^2r_3^2}{A_{123}l_{12}}\kappa_3h_3
$
with $A_{123}=l_{12}l_{13}\sin\theta_1, A_{123}^2=r_1^2r_2^2r_3^2Q.$
	By direct calculations, we have
	\begin{equation}\label{partial h-ijk par kappa3}
		\frac{\partial h_{12,3}}{\partial \kappa_3}=\frac{r_1^2r_2^2r_3^{2}}{A_{123}^3l_{12}}[r_1^2r_2^2r_3^{2}(\kappa_1h_1+\kappa_2h_2)h_3-A_{123}^2(\kappa_1\gamma_2+\kappa_2\gamma_1)].
	\end{equation}
	Note that $v_3$ is the flat vertex of the degenerate triangle $\triangle v_1v_2v_3$ generated by $(r_1,r_2,r_3)$, then $A_{123}=0$ and $h_1>0,h_2>0, h_3<0$ at $(r_1,r_2,r_3)$ by Lemma \ref{basic property I of IDCP} (c), which implies that
	$\frac{\partial h_{12,3}}{\partial \kappa_3}<0$
	around $(r_1,r_2,r_3)$ in the admissible space $\Omega_{123}$ by (\ref{partial h-ijk par kappa3}).  Note that
	$\frac{\partial h_{12,3}}{\partial r_3}=\frac{\partial h_{12,3}}{\partial \kappa_3}\frac{\partial \kappa_3}{\partial r_3}=-\frac{1}{r_3^2}\frac{\partial h_{12,3}}{\partial \kappa_3}.$
	Therefore, there exists $\epsilon>0$ such that $\frac{\partial h_{12,3}}{\partial r_3}(r_1, r_2, r_3+t)>0$
for $t\in (0, \epsilon)$.
	\end{proof}	
\begin{remark}[\cite{Xu MRL} Remark 2.6]\label{limit of h_ij,k}
$h_{ij,k}$ is only defined for non-degenerate inversive distance circle packings in $\Omega_{123}$,
while $h_i$ is defined for any $(r_1, r_2, r_3)\in \mathbb{R}^3_{>0}$. If $(r_1, r_2, r_3)\in \mathbb{R}^3_{>0}$ generates a degenerate triangle $\triangle v_1v_2v_3$
having $v_3$ as the flat vertex, then
	$$h_{12,3}\rightarrow -\infty, h_{13,2}\rightarrow +\infty, h_{23,1}\rightarrow +\infty$$
	as $(\tilde{r}_1, \tilde{r}_2, \tilde{r}_3)\in \Omega_{123}$ tends to $(r_1, r_2, r_3)\in \partial \Omega_{123}$.
If the triangle $\triangle v_1v_2v_3$ generated by a generalized inversive distance circle packing $(r_1, r_2, r_3)$ is degenerate with $v_{k}$ as the flat point,
we denote $h_{ij,k}(r_1, r_2, r_3)=-\infty, h_{ik,j}(r_1, r_2, r_3)=h_{jk,i}(r_1, r_2, r_3)=+\infty$ for simplicity in the following.
By the proof of Lemma \ref{interval} (b),  under the same conditions in Lemma \ref{interval} (b), we further have
\begin{equation*}
\frac{\partial h_{12,3}}{\partial r_3}\rightarrow+\infty, \left(\frac{\partial }{\partial r_3}\left(\frac{h_{12,3}}{l_{12}}\right)\rightarrow+\infty \text{~equivalently}\right)
\end{equation*}
as $(\tilde{r}_1, \tilde{r}_2, \tilde{r}_3)\in \Omega_{123}$ tends to $(r_1, r_2, r_3)\in \partial \Omega_{123}$.
Under the same conditions, one can prove similarly that
\begin{equation*}
  \frac{\partial }{\partial r_3}\left(\frac{h_{13,2}}{l_{13}}\right)\rightarrow-\infty,\frac{\partial }{\partial r_3}\left(\frac{h_{23,1}}{l_{23}}\right)\rightarrow-\infty,
\end{equation*}
as $(\tilde{r}_1, \tilde{r}_2, \tilde{r}_3)\in \Omega_{123}$ tends to $(r_1, r_2, r_3)\in \partial \Omega_{123}$.
\end{remark}

% 	\begin{remark}\label{Vi hi<0, hj,hk>0 Euclidean}
% 		Suppose $\sigma=\{123\}\in F$ is a triangle with a weight $\eta: E_\sigma\rightarrow (-1, +\infty)$
% 		satisfying the structure conditions $(\ref{structure condition})$.
% 		For $(r_1, r_2, r_3)\in V_i$, we have $h_i<0$ and $h_j>0, h_k>0$.
% 	\end{remark}
	
	\subsection{Weighted Delaunay triangulations}
Weighted Delaunay triangulations in Definition \ref{defn of weighted Delaunay} are
natural generalizations of the classical Delaunay triangulations. They have wide applications.
See \cite{AK, Glickenstein DCG, Glickenstein preprint} and others.
In this subsection, we propose an alternative characterization of
weighted Delaunay triangulations for inversive distance circle packing metrics,
and generalize the Definition \ref{defn of weighted Delaunay} of weighted Delaunay triangulations for non-degenerate inversive distance circle packing metrics to generalized inversive distance circle packing metrics.
	
\begin{figure}[!htb]
\centering
  \includegraphics[height=.4\textwidth,width=.6\textwidth]{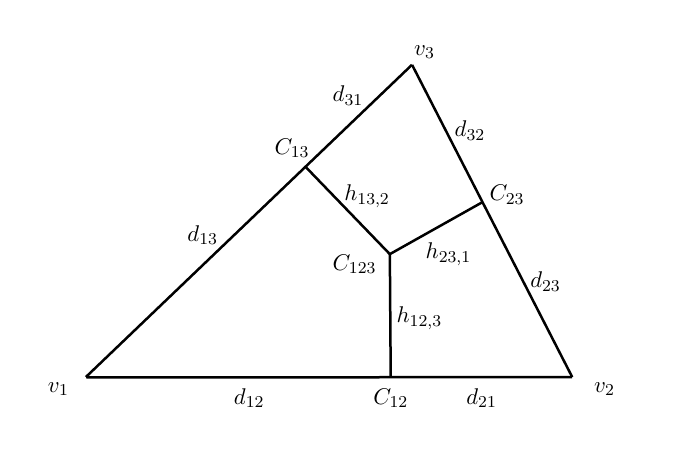}
  \caption{Signed distances of geometric centers}
  \label{Signed distance of geometric center}
\end{figure}

	Assume $r: V\rightarrow (0, +\infty)$ is a non-degenerate inversive distance circle packing on
a weighted triangulated surface $(S, \mathcal{T}, I)$ with the weight $I: E\rightarrow (-1, +\infty)$ satisfying the structure condition (\ref{structure condition}).
Let $\triangle v_1v_2v_3$ be a Euclidean triangle in the plane isometric to a face in $(S, \mathcal{T}, I, r)$.
Then there exists a unique geometric center $C_{123}$ such that its power distances to $v_1,v_2,v_3$ are all the same.
Projections of the geometric center $C_{123}$ to the lines $v_1v_2, v_1v_3, v_2v_3$ give rise to the geometric centers of these edges, which are denoted by $C_{12}, C_{13}, C_{23}$ respectively.
Please refer to Figure \ref{Signed distance of geometric center}.
Denote $d_{ij}$ as the signed distance of $C_{ij}$ to the vertex $v_i$.
Then we have \cite{Glickenstein JDG}
	\begin{equation}\label{d}
		d_{ij}=\frac{r_i^2+r_ir_jI_{ij}}{l_{ij}}, \quad h_{ij,k} = \frac{d_{ik} - d_{ij}\cos \theta_i}{\sin \theta_i},
	\end{equation}
where $\theta_i$ is the inner angle of the triangle $\triangle v_1v_2v_3$ at $v_i$.
Note that $d_{ij}$ could be defined independent of the existence of the geometric center $C_{ijk}$ by (\ref{d}) and
$h_{ij,k}$ is symmetric in the indices $i,j$, while $d_{ij}$ is not.

    \begin{lemma}\label{d>0}
Assume $r: \{v_1,v_2, v_3\}\rightarrow (0, +\infty)$ is a non-degenerate inversive distance circle packing on a weighted triangle $\triangle v_1v_2v_3$ with the weight
$I: E\rightarrow (-1, +\infty)$ satisfying the structure condition (\ref{structure condition}).
\begin{description}
  \item[(a)] If $d_{ij}\leq 0$, then $h_{ij,k}>0$.
  \item[(b)] For any vertex $v_i$ of the triangle $\triangle v_1v_2v_3$, at most one of $d_{ij}$ and $d_{ik}$ is nonpositive.
\end{description}
	\end{lemma}

    \begin{proof}
    If $d_{ij}\leq 0$,  then $I_{ij} \leq -r_i/r_j <0$ by (\ref{d}), which implies $I_{ij}\in (-1, 0)$ by $I_{ij}\in (-1, +\infty)$. As a result, we have $h_k>0$ by the definition of $h_k$ in (\ref{h_i}) and the structure condition (\ref{structure condition}), which further implies $h_{ij,k}>0$ by (\ref{h_ij,k}).

    If we further have $d_{ik}\leq 0$, similar arguments imply $I_{ik}\in (-1, 0)$, which implies $I_{jk}\in (-1, 0)$ by the structure condition $I_{ij}+I_{ik}I_{jk}\geq 0$ and $I_{jk}\in (-1, +\infty)$.
    Without loss of generality, assume that $I_{ij}$ has the largest absolute value among $I_{ij}, I_{jk}, I_{ik}$.
    Then we have $I_{ij}+I_{ik}I_{jk}<0$ by $I_{ij},I_{ik},I_{jk}\in (-1, 0)$, which contradicts the structure condition (\ref{structure condition}). Therefore, at most one of $d_{ij}$ and $d_{ik}$ is nonpositive.
    \end{proof}

\begin{figure}[!htb]
\centering
  \includegraphics[height=.8\textwidth,width=1\textwidth]{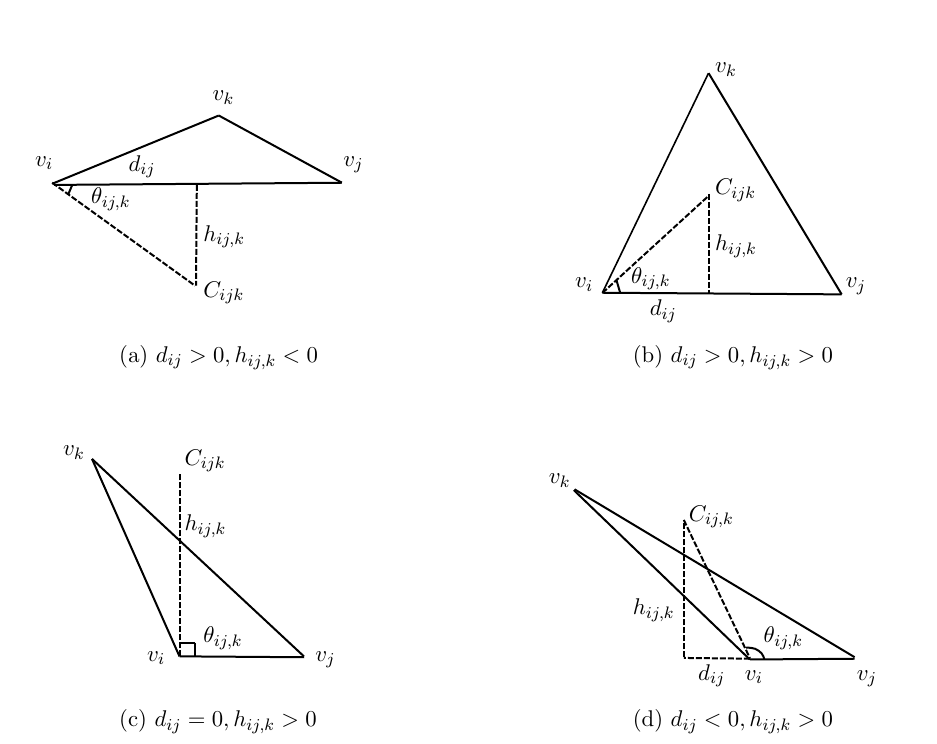}
  \caption{The angle $\theta_{ij,k}$}
  \label{theta ijk}
\end{figure}

\begin{remark}
Lemma \ref{d>0} (b) shows that, for a triangle $\triangle v_1v_2v_3$ generated by a non-degenerate inversive distance circle packing, the geometric center $C_{123}$ can not lie in some regions in the plane determined by $\triangle v_1v_2v_3$.
\end{remark}
Note that weighted Delaunay triangulation in Definition \ref{defn of weighted Delaunay} is only defined
for non-degenerate inversive distance circle packing metrics.
For the following applications, we need to introduce the definition of weighted Delaunay triangulation
for generalized inversive distance circle packing metrics.
To this end, we introduce the following notion.
	\begin{definition}\label{theta}
	Let $r\in \mathbb{R}^V_{>0}$ be a generalized inversive distance circle packing on a weighted triangulated surface $(S, \mathcal{T}, I)$
with the weight $I: E\rightarrow (-1, +\infty)$ satisfying the structure condition (\ref{structure condition}).
$\triangle v_1v_2v_3$ is a generalized triangle in $(S, \mathcal{T}, I, r)$.
If $\triangle v_1v_2v_3$ is non-degenerate, define $\theta_{ij,k}$ as follows
\begin{equation}\label{theta_ij,k for nondege triangle}
  \theta_{ij,k}= \left\{
		\begin{array}{rl}
			\pi+\arctan\frac{h_{ij,k}}{d_{ij}}, & \text{if } d_{ij} < 0,\\
			\frac{\pi}{2}, & \text{if } d_{ij} = 0,\\
			\arctan\frac{h_{ij,k}}{d_{ij}}, & \text{if } d_{ij} > 0.
		\end{array} \right.
\end{equation}
	If $\triangle v_1v_2v_3$ is degenerate,  define $\theta_{ij,k}$ as follows
		\[  \theta_{ij,k}= \left\{
		\begin{array}{rl}
			\frac{\pi}{2}, & \text{if $v_i$ or $v_{j}$ is the flat vertex,}\\
			-\frac{\pi}{2}, & \text{if $v_k$ is the flat vertex.}\\	
		\end{array} \right. \]
	\end{definition}
By definition and Lemma \ref{d>0}, $\theta_{ij,k}\in [-\frac{\pi}{2},\pi)$.
Note that $h_{ij,k}<0$ implies $I_{ij}>1$ by (\ref{h_ij,k}), (\ref{h_i}) and the structure condition (\ref{structure condition}),
which further implies $d_{ij}>0$ by (\ref{d}).
As a result, for a non-degenerate triangle $\triangle v_1v_2v_3$ in $(S, \mathcal{T}, I, r)$,
$\theta_{ij,k}$ is in fact the signed angle $\angle v_jv_iC_{ijk}$
by Lemma \ref{d>0}, which is negative if $h_{ij,k}<0$ and nonnegative otherwise.
Please refer to Figure \ref{theta ijk} for this.
For non-degenerate inversive distance circle packings on a weighted triangle $\triangle v_1v_2v_3$ with
the weight $I: E\rightarrow (-1, +\infty)$ satisfying the structure condition (\ref{structure condition}),
$\theta_{ij,k}$ is obviously a continuous function of $(r_1,r_2,r_3)\in \Omega_{123}$ and satisfies $\theta_{ij,k}+\theta_{ik,j}=\theta_i$ by Lemma \ref{d>0}.
Specially, if  $(r_1,r_2,r_3)\in \Omega_{123}$ tends to $(\bar{r}_1,\bar{r}_2,\bar{r}_3)\in \Omega_{123}$ with
$d_{ij}(\bar{r}_1,\bar{r}_2,\bar{r}_3)=0$, we have $I_{ij}\in (-1,0)$ by (\ref{d}), $h_{ij,k}(\bar{r}_1,\bar{r}_2,\bar{r}_3)>0$
by (\ref{h_ij,k})
and then
$\theta_{ij,k}(r_1,r_2,r_3)\rightarrow \frac{\pi}{2}=\theta_{ij,k}(\bar{r}_1,\bar{r}_2,\bar{r}_3)$ by Definition \ref{theta}.
We further have the following property on $\theta_{ij,k}$ for generalized inversive distance circle packings on a weighted triangle with the weight in $(-1, +\infty)$.

\begin{lemma}\label{theta continuous}
  Suppose $\triangle v_1v_2v_3$ is a face in a weighted triangulated surface $(S, \mathcal{T}, I)$
with the weight $I: E\rightarrow (-1, +\infty)$ satisfying the structure condition (\ref{structure condition}).
Then $\theta_{ij,k}(r_1, r_2, r_3)$ is a continuous function defined on  $\overline{\Omega_{123}}$
and satisfies
\begin{equation}\label{theta++theta-=theta}
\theta_{ij,k}+\theta_{ik,j}=\theta_i.
\end{equation}
\end{lemma}
\begin{proof}
  We just need to prove that $\theta_{ij,k}(r_1,r_2,r_3)\rightarrow \theta_{ij,k}(\bar{r}_1,\bar{r}_2,\bar{r}_3)$
as $(r_1,r_2,r_3)\in \Omega_{123}$ tends to a point $(\bar{r}_1,\bar{r}_2,\bar{r}_3)\in \partial\Omega_{123}$.

If $v_k$ is the flat point of the degenerate triangle $\triangle v_1v_2v_3$ generated by $(\bar{r}_1,\bar{r}_2,\bar{r}_3)$, then $I_{ij}>1$
by Lemma \ref{basic property I of IDCP} (c), which implies $d_{ij}>0$ by (\ref{d}).
Note that $h_{ij,k}(r_1,r_2,r_3)\rightarrow -\infty$ as $(r_1,r_2,r_3)\rightarrow(\bar{r}_1,\bar{r}_2,\bar{r}_3)$ by Remark \ref{limit of h_ij,k}, we have $\theta_{ij,k}(r_1,r_2,r_3)=\arctan\frac{h_{ij,k}}{d_{ij}}\rightarrow -\frac{\pi}{2}=\theta_{ij,k}(\bar{r}_1,\bar{r}_2,\bar{r}_3)$
by Definition \ref{theta}.

If $v_i$ is the flat point of the degenerate triangle $\triangle v_1v_2v_3$ generated by $(\bar{r}_1,\bar{r}_2,\bar{r}_3)$,
then $h_{ij,k}(r_1,r_2,r_3)\rightarrow +\infty$ as $(r_1,r_2,r_3)\rightarrow(\bar{r}_1,\bar{r}_2,\bar{r}_3)$ by Remark \ref{limit of h_ij,k}.
As a result, we have $\theta_{ij,k}(r_1,r_2,r_3)\rightarrow \frac{\pi}{2}=\theta_{ij,k}(\bar{r}_1,\bar{r}_2,\bar{r}_3)$
as $(r_1,r_2,r_3)\rightarrow(\bar{r}_1,\bar{r}_2,\bar{r}_3)$
by Definition \ref{theta},
no matter the sign of $d_{ij}(r_1,r_2,r_3)$.
The same argument applies to the case that $v_j$ is the flat point.
\end{proof}

The definition of weighted Delaunay triangulation for inversive distance circle packings has the following relationships with $\theta_{ij,k}$.
	\begin{lemma}\label{g-Delaunay lemma}
Suppose $r\in \mathbb{R}^V_{>0}$ is a non-degenerate inversive distance circle packing on a weighted triangulated surface $(S, \mathcal{T}, I)$
with the weight $I: E\rightarrow (-1, +\infty)$ satisfying the structure condition (\ref{structure condition}).
An edge $v_iv_j\in E$ is shared by two adjacent non-degenarate triangles $\triangle v_iv_jv_k$
and $\triangle v_iv_jv_l$ in $(S, \mathcal{T}, I, r)$.
Then the edge $v_iv_j$ is weighted Delaunay in the inversive distance circle packing metric if and only if
$$\theta_{ij,k}+\theta_{ij,l}\ge0.$$
Furthermore, if $d_{ij}\le0$, then $h_{ij,k}>0, h_{ij,l}>0$, $\theta_{ij,k}\geq \frac{\pi}{2},\theta_{ij,l}\geq \frac{\pi}{2}$,
which implies $h_{ij,k}+h_{ij,l}>0$ and $\theta_{ij,k}+\theta_{ij,l}\geq \pi>0$.
	\end{lemma}
    \begin{proof}
	If $d_{ij}>0$, then $\theta_{ij,k}=\arctan\frac{h_{ij,k}}{d_{ij}}\in (-\frac{\pi}{2},\frac{\pi}{2})$ and
$\theta_{ij,l}=\arctan\frac{h_{ij,l}}{d_{ij}}\in (-\frac{\pi}{2},\frac{\pi}{2})$ by Definition \ref{theta}. In this case, we have $$\frac{h_{ij,k}+h_{ij,l}}{d_{ij}}=\tan\theta_{ij,k}+\tan\theta_{ij,l}=\frac{\sin(\theta_{ij,k}+\theta_{ij,l})}{\cos\theta_{ij,k}\cos\theta_{ij,l}},$$
which implies	$h_{ij,k}+h_{ij,l}\ge0$ is equivalent to $\theta_{ij,k}+\theta_{ij,l}\ge0$.
If $d_{ij}\le0$, we have $h_{ij,k}>0$ and $h_{ij,l}>0$ by Lemma \ref{d>0}, which implies $\theta_{ij,k}\ge\frac{\pi}{2}>0$ and $\theta_{ij,l}\ge\frac{\pi}{2}>0$ by Definition \ref{theta}. Therefore, $h_{ij,k}+h_{ij,l}>0$ and $\theta_{ij,k}+\theta_{ij,l}\geq \pi>0$.
    \end{proof}

Note that $h_{ij,k}$ is only defined for non-degenerate inversive distance circle packing metrics,
while $\theta_{ij,k}$ could be defined for generalized inversive distance circle packing metrics.
Motivated by Lemma \ref{g-Delaunay lemma}, we introduce the following definition of weighed Delaunay triangulation
for generalized inversive distance circle packing metrics, which generalizes Definition \ref{defn of weighted Delaunay} of
weighed Delaunay triangulation
for non-degenerate inversive distance circle packing metrics.

\begin{definition}\label{defn of weighted Delaunay for generalized IDCP}
Suppose $r: V\rightarrow (0, +\infty)$ is a generalized inversive distance circle packing
on a weighted triangulated surface $(S, \mathcal{T}, I)$
with the weight $I: E\rightarrow (-1, +\infty)$ satisfying the structure condition (\ref{structure condition}).
$v_iv_j\in E$ is an interior edge shared by two adjacent triangles $\triangle v_iv_jv_k$ and $\triangle v_iv_jv_m$ in $\mathcal{T}$.
$v_iv_j\in E$ is weighted Delaunay in the generalized inversive distance circle packing $r$ on $(S, \mathcal{T}, I)$ if $$\theta_{ij,k}+\theta_{ij,m}\ge0.$$
The triangulation $\mathcal{T}$ is weighted Delaunay in the generalized inversive distance circle packing $r$ on $(S, \mathcal{T}, I)$ if every interior edge is weighed Delaunay in $r$.
\end{definition}

For simplicity, we also say $r$
is a generalized weighted Delaunay inversive distance circle packing
on $(S, \mathcal{T}, I)$, if $\mathcal{T}$ is weighted Delaunay in $r$.

    \begin{lemma}\label{monotonicity}
Suppose $\triangle v_1v_2v_3$ is a face in a weighted triangulated surface $(S, \mathcal{T}, I)$
with the weight $I: E\rightarrow (-1, +\infty)$ satisfying the structure condition (\ref{structure condition}),
and $ (r_1, r_2, \hat{r}_3)$ and $ (r_1, r_2, \bar{r}_3)$ are two generalized inversive distance circle packings on $\triangle v_1v_2v_3$ with $\hat{r}_3<\bar{r}_3$. If $r_1,r_2$ are fixed, $d_{12}(r_1, r_2)>0$ and $\Delta_{123}>0$, then $\theta_{12,3}$ is strictly increasing in $r_3\in [\hat{r}_3,\bar{r}_3]$.
    \end{lemma}
    \begin{proof}
    By Lemma \ref{interval} (a), $ (r_1, r_2, r_3)$ generates a non-degenerate triangle $\triangle v_1v_2v_3$ for $r_3\in (\hat{r}_3,\bar{r}_3)$.
For $r_3\in (\hat{r}_3,\bar{r}_3)$, $h_{12,3}$ and $\theta_{12,3}$ are smooth functions of $r_3$.
By direct calculations, we have
    \begin{equation}\label{hi-}
     \begin{aligned}
\frac{\partial{h_{12,3}}}{\partial{\kappa_3}}&=\frac{r_1^2r_2^2r_3^2}{A_{123}^3l_{12}}[r_1^2r_2^2r_3^2(\kappa_1h_1+\kappa_2h_2)h_3-A_{123}^2(\kappa_2\gamma_1+\kappa_1\gamma_2)] \notag\\
     	&=\frac{ r_1^4r_2^4r_3^3}{A_{123}^3l_{12}}(1-I_{12}^2-I_{13}^2-I_{23}^2-2I_{12}I_{13}I_{23})(\kappa_1^2+\kappa_2^2+2\kappa_1\kappa_2I_{12}) \\
     	& = -\frac{ r_1^4r_2^4r_3^3}{A_{123}^3l_{12}}\Delta_{123}(\kappa_1^2+\kappa_2^2+2\kappa_1\kappa_2I_{12})\\
       &<0
     \end{aligned}
     \end{equation}
by $\Delta_{123}>0$ and $I_{12}>-1$,
which further implies
    $$\frac{\partial{\theta_{12,3}}}{\partial r_3}=-\frac{d_{12}\kappa_3^{2}}{d_{12}^2+(h_{12,3})^2}\cdot\frac{\partial{h_{12,3}}}{\partial{\kappa_3}}>0,\ \ \forall r_3\in (\hat{r}_3,\bar{r}_3)$$
by the definition of $\theta_{12,3}$ and the assumption $d_{12}=d_{12}(r_1, r_2)>0$.  Note that $\theta_{12,3}$ is a continuous function of $r_3\in [\hat{r}_3,\bar{r}_3]$ by Lemma \ref{theta continuous},
we have $\theta_{12,3}$ is strictly increasing in $r_3\in [\hat{r}_3,\bar{r}_3]$.
    \end{proof}

\section{Maximal principle, Ring Lemma and spiral hexagonal triangulations}\label{section 3}
In this section, we derive a maximal principle for generic inversive distance circle packings,
which is a generalization of the maximal principle obtained in \cite{H2} for Thurston's circle packings. Then we give a ring lemma for inversive distance circle packings in the hexagonal triangulated plane with inversive distance $I: E\rightarrow (-\frac{1}{2}, +\infty)$, which generalizes the ring lemma obtained for Thurston's circle packings in \cite{H2} in the hexagonal triangulated plane. We further obtain some properties of the linear discrete conformal factors of inversive distance circle packings on the hexagonal triangulated plane.
\subsection{Maximal principle}
Let $P_n$ be a star-shaped $n$-sided polygon in the plane with boundary vertices $v_1,\cdots,v_n$ cyclically ordered ($v_{n+i} = v_i$).
Assume $v_0$ is an interior point of $P_n$ and it induces a triangulation $\mathcal{T}$ of $P_n$ with triangles $ v_0v_iv_{i+1}$.
We take the assignment of radii $r:V(\mathcal{T}) \to \mathbb{R}_{>0}$ as a vector in $\mathbb{R}^{n+1}$.
 For any two vectors $x = (x_0,\dots, x_n)$ and $y = (y_0,\dots, y_n)$ in $\mathbb{R}^{n+1}$,
 we use $x\geq y$ to denote $x_i\geq y_i$ for all $i \in \{0, \dots, n\}$.

   \begin{figure}[!htb]
\centering
  \includegraphics[height=.4\textwidth,width=.5\textwidth]{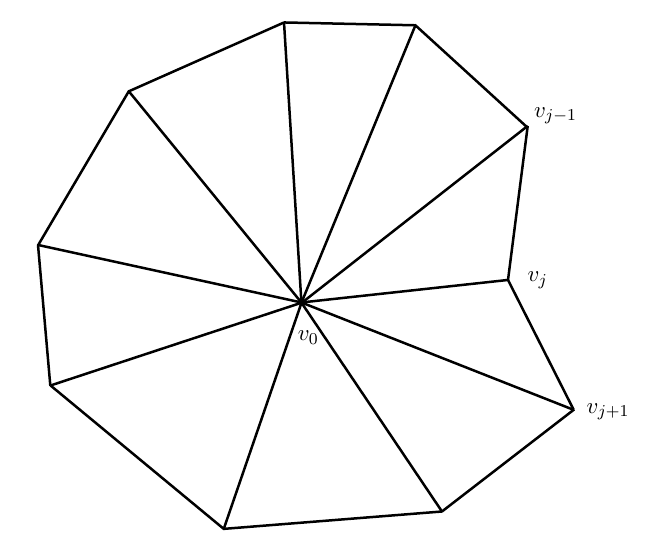}
  \caption{star triangulation of a polygon}
  \label{star triangulation}
\end{figure}

We have the following maximal principle for generic inversive distance circle packings.

\begin{theorem}[Maximal principle]\label{Maximum principle}
  Let $\mathcal{T}$ be a star triangulation of $P_n$ with boundary vertices $v_1,\dots, v_n$ and a unique interior vertex $v_0$.
  $I$ is a regular weight defined on the edges in $\mathcal{T}$ satisfying the structure condition (\ref{structure condition})
  with $I:E\rightarrow (-1, 1]$ or $I:E\rightarrow [0, +\infty)$.
  If $\overline{r}$ and $r$ are two generalized  inversive distance circle packings on $(P_n, \mathcal{T}, I)$ satisfying
  \begin{description}
    \item[(a)] $\overline{r}$ and $r$ are  generalized weighted Delaunay inversive distance circle packings on $(P_n, \mathcal{T}, I)$,
    \item[(b)] the combinatorial curvatures $K_0(r)$ and $K_0(\bar{r})$ at the vertex $v_0$ satisfy $K_0(r)\le K_0(\bar{r})$,
    \item[(c)] $\max\{\frac{r_i}{\bar{r}_i}|i=1,2,\dots,n\}\le\frac{r_0}{\bar{r}_0}$,
  \end{description}
  then $\frac{r_i}{\bar{r}_i}=const$ for any $i=0, 1,\dots,n$.
\end{theorem}

 We will use the following notations to prove Theorem \ref{Maximum principle}.
 For $i\in \{1,\cdots, n\}$, we denote $I_{0i}$ as $I_i$ for simplicity.
 For two adjacent triangles $ \triangle v_0v_jv_{j\pm1}$ in $\mathcal{T}$, set
	$$\Delta_i^-=\Delta_{0i(i-1)}, \Delta_i^+=\Delta_{0i(i+1)},
	h_j^-=h_{0j,j-1}, h_j^+=h_{0j,j+1},
	\theta_j^-=\theta_{0j,j-1}, \theta_j^+=\theta_{0j,j+1}.$$

The proof of maximal principle is based on the following key lemma.

  \begin{lemma}\label{key lemma}
  	If $r, \overline{r}:\{v_0,v_1,\dots,v_n\}\rightarrow\mathbb{R}_{>0}$ satisfy (a), (b), (c) in Theorem \ref{Maximum principle} and there exists $j\in\{1,2,\dots,n\}$ such that $\frac{r_j}{\bar{r}_j}<\frac{r_0}{\bar{r}_0}$, then there exists $\hat{r}\in\mathbb{R}_{>0}^{n+1}$ such that
  \begin{description}
    \item[(a)] $\hat{r}_i\geq r_i$ for $i\in \{1,\cdots, n\}$,
    \item[(b)] $\frac{\hat{r}_i}{\bar{r}_i}\le\frac{\hat{r}_0}{\bar{r}_0}=\frac{r_0}{\bar{r}_0}$ for all $i=1,2,\dots,n$,
    \item[(c)] $\hat{r}$ is a generalized weighted Delaunay inversive distance circle packing on $(P_n, \mathcal{T}, I)$,
    \item[(d)] let $\alpha(r)$ be the cone angle of the inversive distance circle packing $r$ at $v_0$, then
  \begin{equation}
  \alpha(\hat{r})>\alpha(r).
  \end{equation}
  \end{description}
  \end{lemma}
  \begin{proof}
Without loss of generality, we may assume that $r_0=\bar{r}_0$, otherwise we can scale
$r_i$ ($i\in \{0, \cdots, n\}$) by a factor $\frac{\overline{r}_0}{r_0}$.
Then the condition (c) in Theorem \ref{Maximum principle} is equivalent to $r_i\le{\bar{r}_i}$ for all $i\in\{1,2,\dots,n\}$.

If the weight $I$ takes all the value in $(-1, 1]$, i.e. $I: E\rightarrow (-1, 1]$, then the triangles $\triangle v_0v_iv_{i+1}$, $i\in \{1, \cdots, n\}$, are non-degenerate
for any $r\in \mathbb{R}^{n+1}_{>0}$ by Lemma \ref{basic property I of IDCP} (a).
Furthermore, we have $h_i^+(r)\geq 0$ and $h_i^-(r)\geq 0$ for any $i\in \{1, \cdots, n\}$ and $r\in \mathbb{R}^{n+1}_{>0}$
by $I: E\rightarrow (-1, 1]$ and Lemma \ref{basic property I of IDCP} (b), which implies
$h_i^+(r)+h_i^-(r)\geq 0$.
If $h_i^+(r)+h_i^-(r)=0$, then $h_i^+(r)=0$ and $h_i^-(r)=0$, which implies that $I_{0i}=1$, $I_{0,i+1}=-I_{i, i+1}\in (-1, 1)$
and $I_{0,i-1}=-I_{i, i-1}\in (-1, 1)$ by Lemma \ref{basic property I of IDCP} (b). This contradicts the assumption that the weight $I$ is regular.
Therefore, $h_i^+(r)+h_i^-(r)>0$ for any $i\in \{1, \cdots, n\}$ and $r\in \mathbb{R}^n_{>0}$.
Specially, $h_j^+(r)+h_j^-(r)>0$ for $j\in \{1, \cdots, n\}$ with $r_j<\bar{r}_j$.
As a result, we have $\frac{\partial \alpha}{\partial r_j}>0$ by Lemma \ref{basic property I of IDCP} (b)
and then $\alpha(\hat{r})>\alpha(r)$ for $\hat{r}=(r_1, \cdots, r_{j-1}, r_j+t, r_{j+1}, \cdots, r_n), t\in (0, \bar{r}_j-r_j)$.

If the weight $I$ takes all the value in $[0, +\infty)$, i.e. $I: E\rightarrow [0, +\infty)$, set
\begin{equation*}
  \begin{aligned}
J=&\{j\in\{1,2,\dots,n\}|r_j<{\bar{r}_j}\}, \\
K=&\{k\in\{1,2,\dots,n\}|r_k={\bar{r}_k}\},\\
\gamma(r)=&\sum_{j\in{J}}(\theta_{0j,j+1}+\theta_{0j,j-1})=\sum_{j\in{J}}(\theta_j^++\theta_j^-), \\
\beta(r)=&\sum_{k\in{K}}(\theta_{0k,k+1}+\theta_{0k,k-1})=\sum_{k\in{K}}(\theta_k^++\theta_k^-).
  \end{aligned}
\end{equation*}
Then $J \neq \emptyset$ by assumption. By (\ref{theta++theta-=theta}), we have
$
\alpha(r)=\beta(r)+\gamma(r), \alpha(\overline{r})=\beta(\overline{r})+\gamma(\overline{r}),
$
which further implies
		\begin{equation}\label{beta+gamma<beta+gamma}
			\beta(\bar{r})+\gamma(\bar{r})\leq \beta(r)+\gamma(r).
		\end{equation}
by the condition $K_0(r)\le K_0(\bar{r})$.

        \textbf{Claim 1}: For any $j\in J$, $\theta_{0}^{j,j-1}(r)<\pi$ and $\theta_{0}^{j,j+1}(r)<\pi$.

        We just need to prove that for any $j\in J$, $v_0$ is not the flat point if the triangle $\triangle v_0v_{j}v_{j-1}$ generated by $r$ is degenerate.
        Otherwise, suppose for some $j\in J$, $v_0$ is the flat vertex of the degenerate triangle $\triangle v_0v_{j}v_{j-1}$ generated by $r$.
        Then $I_{j,j-1}>1$ by Lemma \ref{basic property I of IDCP} (c), which further implies that $\Delta_{0,j-1,j}=I_j^2+I_{j-1}^2+I_{j,j-1}^2+2I_jI_{j-1}I_{j,j-1}-1>0$ by Lemma \ref{basic property I of IDCP} (d).
        By Lemma \ref{simply connect of admi space with weight}, $r$ satisfies $\kappa_0=f(\kappa_{j-1},\kappa_{j})$, where
        \begin{equation*}
        \begin{aligned}
        f(\kappa_{j-1},\kappa_{j})
        =\frac{1}{I_{j,j-1}^2-1}\left[(\kappa_j\gamma_{j-1}+\kappa_{j-1}\gamma_j)+
        \sqrt{\Delta_{0,j-1,j}(\kappa_{j}^2+\kappa_{j-1}^2+2\kappa_j\kappa_{j-1}I_{j,j-1})}\right]
        \end{aligned}
        \end{equation*}
with $\gamma_j=I_{0,j-1}+I_{0,j}I_{j,j-1}\geq 0$ and $\gamma_{j-1}=I_{0,j}+I_{0,j-1}I_{j,j-1}\geq 0$ by the structure condition (\ref{structure condition}).
Note that $\kappa_j>\overline{\kappa}_j$ and $\kappa_{j-1}\geq\overline{\kappa}_{j-1}$, we have
$$\overline{\kappa}_0=\kappa_0=f(\kappa_{j-1},\kappa_{j})>f(\overline{\kappa}_{j-1},\overline{\kappa}_{j}).$$
This implies that $(\overline{r}_0, \overline{r}_j, \overline{r}_{j-1})$ is in the complement of the space of generalized inversive distance circle packings on $\triangle v_0v_{j}v_{j-1}$ in $\mathbb{R}^3_{>0}$ by Lemma \ref{simply connect of admi space with weight}, which contradicts the assumption that $\overline{r}$ is a generalized inversive distance circle packing on $(P_n, \mathcal{T}, I)$.

		\textbf{Claim 2}: There exists $j\in J$ such that $\theta_j^+(r)+\theta_j^-(r)>0$.

%		To prove the claim, we just need to consider the case that $d_{0j}(r)>0$ for all $j\in J$.
%Otherwise, there exists some $j\in J$ such that $d_{0j}(r)\leq 0$, which implies $I_j\in (-1,0)$ by (\ref{d}).
%As a result, if the triangle $\triangle v_0v_jv_{j\pm1}$ is non-degenerate, we have $h_{j}^{\pm}(r)>0$ by Lemma \ref{d>0},
%which implies $\theta_j^\pm(r)\geq \frac{\pi}{2}$ by Definition \ref{theta}.
%If the triangle $\triangle v_0v_jv_{j-1}$ degenerates, then $v_0$ is not the flat vertex by Claim 1 and $v_{j-1}$ is not the flat vertex by
%$I_j\in (-1,0)$ and Lemma \ref{basic property I of IDCP} (c), which implies that $v_j$ is the flat vertex of the degenerate triangle $\triangle v_0v_jv_{j-1}$.
%In this case, we have $h_{j}^-(r)=+\infty$ by Remark \ref{limit of h_ij,k}, which implies $\theta_j^-(r)=\frac{\pi}{2}$ by Definition \ref{theta}.
%Similar arguments show that $\theta_j^+(r)=\frac{\pi}{2}$ if the triangle $\triangle v_0v_jv_{j+1}$ degenerates.
%Therefore, $\theta_j^+(r)+\theta_j^-(r)\geq \pi>0$ if $d_{0j}(r)\leq 0$.

To prove Claim 2, we just need to consider the cases $K\ne\emptyset$ and $K=\emptyset$.
\begin{description}
  \item[Case 1] $K\ne\emptyset$.
\end{description}

		If $K\ne\emptyset$, there exists $i\in K$ such that $i-1$ or $i+1$ is in $J$ as $J\neq \emptyset$.
We just need to consider the following two subcases.
\begin{description}
  \item[Case 1(a)] for any $i\in K$, we have $\Delta_i^->0$ when $r_{i-1}<\bar{r}_{i-1}$, and $\Delta_i^+>0$ when $r_{i+1}<\bar{r}_{i+1}$.
  \item[Case 1(b)] there exists a vertex $i\in K$ such that $\Delta_i^-\le0$ with $r_{i-1}<\bar{r}_{i-1}$ or $\Delta_i^+\le0$ with $r_{i+1}<\bar{r}_{i+1}$.
\end{description}

		In Case 1(a), we have $d_{0i}(r)>0$ by $I\in [0, +\infty)$.  By Lemma \ref{monotonicity}, for any $i\in K$,  $\theta_i^-$ and $\theta_i^+$ are strictly increasing in $r_{i-1}$ and $r_{i+1}$ respectively,
which implies that $\beta(r)\leq\beta(\bar{r})$.
As $J\neq \emptyset$, there exists $i\in K$ such that $i-1$ or $i+1$ is in $J$. Say $i-1\in J$, then $r_{i-1}<\bar{r}_{i-1}$ and then $\theta_i^-(r)<\theta_i^-(\overline{r})$ by Lemma \ref{monotonicity}.
Thus, $\beta(r)<\beta(\bar{r})$, which implies $\gamma(r)>\gamma(\bar{r})\geq 0$ by (\ref{beta+gamma<beta+gamma}).
		Therefore, there exists $j\in J$ such that $\theta_j^+(r)+\theta_j^-(r)>0$ by the definition of $\gamma(r)$.

		In Case 1(b), without loss of generality, we assume that there exists $i_0\in K$, $i_0-1\in J$ such that $\Delta_{i_0}^-\le0$ with $r_{i_0-1}<\bar{r}_{i_0-1}$. By Lemma \ref{basic property I of IDCP} (d), for the triangle $\triangle v_0v_{i_0}v_{i_0-1}$, we have $I_{i_0-1}\in[0,1]$ and the triangle $\triangle v_0v_{i_0}v_{i_0-1}$ is non-degenerate for any $r\in \mathbb{R}^{n+1}_{>0}$.

If $I_{i_0-1}\in[0,1)$, we have $h_{i_{0}-1}^{+}(r)>0$ by (\ref{h_ij,k}) and the structure condition (\ref{structure condition}).
For the triangle $\triangle v_0v_{i_0-1}v_{i_0-2}$, it is non-degenerate or degenerate with $v_{i_0-1}$ as the flat vertex by Claim 1 and
Lemma \ref{basic property I of IDCP} (c), in which cases we have $h_{i_{0}-1}^{-}(r)>0$ by (\ref{h_ij,k}) and the structure condition (\ref{structure condition})
and $h_{i_{0}-1}^{-}(r)=+\infty$ by Remark \ref{limit of h_ij,k} respectively.
Note that $d_{0,i_0-1}(r)>0$, we have $\theta_{i_0-1}^{\pm}(r)>0$ by Definition \ref{theta}, which implies
$\theta_{i_0-1}^{+}(r)+\theta_{i_0-1}^{-}(r)>0$.

If $I_{i_0-1}=1$, by $\Delta_{i_0}^-\le0$, we have
		\begin{equation*}
\begin{aligned}
	0\geq I_{i_0,i_0-1}^2+I_{i_0-1}^2+I_{i_0}^2+2I_{i_0,i_0-1}I_{i_0-1}I_{i_0}-1
=(I_{i_0,i_0-1}+I_{i_0})^2
\geq 0,
\end{aligned}
		\end{equation*}
which implies $I_{i_0,i_0-1}=-I_{i_0}\in [0,1)$ and then $I_{i_0,i_0-1}=-I_{i_0}=0$.
Then $h_{i_0-1}^+(r)=0$ by (\ref{h_ij,k}) and $\theta_{i_0-1}^+(r)=0$ by $d_{0,i_0-1}(r)>0$.
As $I_{i_0-1}=1$, the triangle $\triangle v_0v_{i_0-1}v_{i_0-2}$ is non-degenerate or degenerate with $v_{i_0-1}$ as the flat vertex by
Claim 1 and Lemma \ref{basic property I of IDCP} (c), in which cases $h_{i_0-1}^-(r)\geq0$ by (\ref{h_ij,k}) and the structure condition (\ref{structure condition}) and $h_{i_0-1}^-(r)=+\infty$ by Remark \ref{limit of h_ij,k} respectively.
If $h_{i_0-1}^-(r)=0$, we have $I_{i_0-2}=-I_{i_0-1,i_0-2}=0$ by (\ref{h_ij,k}) and the structure condition (\ref{structure condition}),
which contradicts the assumption that the weight $I$ is regular. Therefore, $h_{i_0-1}^-(r)>0$ or $h_{i_0-1}^-(r)=+\infty$, which implies
$\theta_{i_0-1}^-(r)>0$ by $d_{0,i_0-1}(r)>0$ and Definition \ref{theta}. Therefore, $\theta_{i_0-1}^{+}(r)+\theta_{i_0-1}^{-}(r)=\theta_{i_0-1}^{-}(r)>0$ in this case.
\begin{description}
  \item[Case 2] $K=\emptyset$.
\end{description}

		If $K=\emptyset$, we have $J=\{1,\dots,n\}$,
		\begin{equation}
			\gamma(r)=\sum_{j\in{J}}(\theta_j^+(r)+\theta_j^-(r))=\alpha(r)\ge0.
		\end{equation}

		If $\alpha(r)>0$, there exists $j\in J$ such that $\theta_j^+(r)+\theta_j^-(r)>0$.

		If $\alpha(r)=0$, for any triangle $\triangle v_0v_jv_{j-1}$, $j=1,\dots,n$, the inner angle at vertex $v_0$ is equal to $0$. Thus all triangles are degenerate and flat vertices are not $v_0$. We rule out the case that $I_{j}>1$ for all $j\in J=\{1,\dots,n\}$.
Otherwise, for any triangle $\triangle v_0v_jv_{j-1}$, the flat vertex is $v_j$ or $v_{j-1}$ by Claim 1. Then $\{\theta_j^-(r), \theta_{j-1}^+(r)\}=\{\frac{\pi}{2}, -\frac{\pi}{2}\}, \forall j\in \{1, \cdots, n\}.$
		Without loss of generality, we may assume $v_1$ is the flat vertex of triangle $\triangle v_0v_1v_2$. Then $\theta_1^+(r)=\frac{\pi}{2}$, $\theta_2^-(r)=-\frac{\pi}{2}$ by Definition \ref{theta} and $l_{02}(r)=l_{01}(r)+l_{12}(r)>l_{01}(r).$
		By the weighted Delaunay condition (a) in Theorem \ref{Maximum principle}, $\theta_2^+(r)=\frac{\pi}{2}$, which implies $\theta_3^-(r)=-\frac{\pi}{2}$ and
		$l_{03}(r)=l_{02}(r)+l_{23}(r)>l_{02}(r).$
		By induction, we have
		$$l_{01}(r)<l_{02}(r)<\dots<l_{0n}(r)<l_{01}(r),$$
		which is impossible. So there exists $j\in J$ such that $I_{j}\in[0,1]$.  By Claim 1 and Lemma \ref{basic property I of IDCP} (c),
the flat vertex of the degenerate triangles $\triangle v_0v_jv_{j\pm1}$ is $v_{j}$, which implies $h_j^\pm(r)=+\infty$ by Remark \ref{limit of h_ij,k} and $\theta_{j}^+(r)=\theta_{j}^-(r)=\frac{\pi}{2}$ by Definition \ref{theta}. Therefore, $\theta_{j}^+(r)+\theta_{j}^-(r)=\pi>0$.
This completes the proof of Claim 2.

	Now we fix $j\in J$ in Claim 2. Then we have
		\begin{equation}\label{g-Delaunay}
			\theta_j^+(r)+\theta_j^-(r)>0.
		\end{equation}
		In the following, we will show that there exists $\epsilon>0$
such that $\hat{r}=(r_0,\dots,r_j+t,\dots,r_n)$ satisfies Lemma \ref{key lemma}
for $t\in (0, \epsilon)$.
It is easy to check that for $t\in (0, \overline{r}_j-r_j)$, $\hat{r}$ satisfies Lemma \ref{key lemma} (a) and (b).

		To see part (c) of Lemma \ref{key lemma}, we first show that there exists $\epsilon>0$ such that
$\hat{r}$ is a generalized inversive distance circle packing on $(P_n, \mathcal{T}, I)$ for $t\in (0,\epsilon)$.
Furthermore, we will show that the triangles $\triangle v_0v_jv_{j\pm1}$ generated by $\hat{r}$ are non-degenerate.

%If $d_{0j}(r)\leq 0$, then we have $I_{j}\in (-1,0)$ by (\ref{d}), which implies the triangles $\triangle v_0v_jv_{j\pm 1}$ generated by $r$ are non-degenerete
%or degenerate with $v_j$ as the flat vertex by Claim 1 and Lemma \ref{basic property I of IDCP} (c). By Lemma \ref{interval} (b),
%there exists $\epsilon>0$ such that for $t\in (0, \epsilon)$, $\hat{r}$ is a generalized inversive distance circle packing on $(P_n, \mathcal{T}, I)$
%and the triangles $\triangle v_0v_jv_{j\pm1}$ generated by $\hat{r}$ are non-degenerate.

The triangle $\triangle v_0v_jv_{j-1}$ generated by $r$ is non-degenerate or degenerate with $v_j$ or $v_{j-1}$ as the flat vertex by Claim 1.
By Lemma \ref{interval} (b), we just need to prove that $v_{j-1}$ is not the flat vertex of the triangle $\triangle v_0v_jv_{j-1}$ generated by $r$ if it degenerates.
Otherwise, we have $\theta_j^-(r)=-\frac{\pi}{2}$ by Definition \ref{theta}, which implies $\theta_j^+(r)>\frac{\pi}{2}$ by (\ref{g-Delaunay}).
Note that $d_{0j}(r)>0$, we have $\theta_j^+(r)\in [-\frac{\pi}{2}, \frac{\pi}{2}]$ by Definition \ref{theta}, which contradicts $\theta_j^+(r)>\frac{\pi}{2}$.
Therefore, $v_{j-1}$ can never be the flat vertex of the triangle $\triangle v_0v_jv_{j-1}$ if it is degenerate.
Similar arguments applying to the triangle $\triangle v_0v_jv_{j+1}$ show that $v_{j+1}$ can never be
the flat vertex of the triangle $\triangle v_0v_jv_{j+1}$ if it is degenerate.
Therefore, by Lemma \ref{interval} (b), there exists $\epsilon>0$ such that for $t\in (0, \epsilon)$, $\hat{r}=(r_0,\dots,r_j+t,\dots,r_n)$ is a generalized inversive distance circle packing on $(P_n, \mathcal{T}, I)$ and the triangles $\triangle v_0v_jv_{j\pm1}$ generated by $\hat{r}$ are non-degenerate.
		
		Next, we show $\hat{r}$ satisfies the weighted Delaunay condition. As $\hat{r}$ differs from $r$ only at the $j$-th position,
we just need to consider the edges $v_0v_j$ and $v_0v_{j\pm1}$.
For the edge $v_0v_j$, since $\theta_j^+(r)+\theta_j^-(r)>0$, we have $\theta_j^+(\hat{r})+\theta_j^-(\hat{r})>0$ for small $t>0$ by the continuity of $\theta_j^{\pm}$ in Lemma \ref{theta continuous}.
For the edge $v_0v_{j-1}$, $\theta_{j-1}^-(r)=\theta_{j-1}^-(\hat{r})$.
%If $d_{0,j-1}(r)\leq0$, we have $I_{j-1}\in (-1,0)$ by (\ref{d}), which implies that $\theta_{j-1}^{+}(r)\geq \frac{\pi}{2}$ if the triangle $\triangle v_0v_jv_{j-1}$ is non-degenerate by Lemma \ref{g-Delaunay lemma} and $\theta_{j-1}^+(r)=\frac{\pi}{2}$ if the triangle $\triangle v_0v_jv_{j-1}$ is degenerate with $v_{j-1}$ as the flat vertex by Claim 1, Lemma \ref{basic property I of IDCP} (c) and Definition \ref{theta}. The same argument applying to the triangle $\triangle v_0v_{j-1}v_{j-2}$ gives  $\theta_{j-1}^-(r)\geq \frac{\pi}{2}$. Therefore, if $d_{0,j-1}(r)\leq0$, we have
%$\theta_{j-1}^+(r)+\theta_{j-1}^-(r)\geq \pi$. The conclusion then follows from the continuity of $\theta_{j-1}^{\pm}$ in Lemma \ref{theta continuous}.
If $\Delta_j^->0$, we have $\theta_{j-1}^+(r)< \theta_{j-1}^+(\hat{r})$ for $t\in (0, \overline{r}_j-r_j)$ by Lemma \ref{monotonicity}, which implies
		$\theta_{j-1}^+(\hat{r})+\theta_{j-1}^-(\hat{r})>\theta_{j-1}^+(r)+\theta_{j-1}^-(r)\ge0.$
		This implies that the edge $v_0v_{j-1}$ satisfies the weighted Delaunay condition for $\hat{r}$.
If $\Delta_{j}^-\le0$, we have $I_{j-1}\in [0,1]$ and the triangle $\triangle v_0v_jv_{j-1}$ generated by any $r\in \mathbb{R}^{n+1}_{>0}$ is non-degenerate by Lemma \ref{basic property I of IDCP} (c) (d).
Repeat the arguments in Case 1 (b) in the proof of Claim 2, we have $\theta_{j-1}^+(r)+\theta_{j-1}^-(r)>0$ if $I_{j-1}\in [0,1)$,
and $\theta_{j-1}^+(\hat{r})+\theta_{j-1}^-(\hat{r})=\theta_{j-1}^+(r)+\theta_{j-1}^-(r)=\theta_{j-1}^-(r)>0$ if $I_{j-1}=1$.
The conclusion then follows from the continuity of $\theta_{j-1}^{\pm}$ in Lemma \ref{theta continuous}.
Therefore, there exists $\epsilon>0$ such that the edge $v_0v_{j-1}$ satisfies the weighed Delaunay condition in $\hat{r}$ for $t\in (0,\epsilon)$.
The same arguments apply to the edge $v_0v_{j+1}$.
		
		To see part (d) of Lemma \ref{key lemma}, by the arguments for part (c), there exists $\epsilon>0$ such that
the triangles $\triangle v_0v_jv_{j\pm1}$ are non-degenerate in $\hat{r}$ and $\theta_{j}^+(\hat{r})+\theta_{j}^-(\hat{r})>0$
for $t\in (0,\epsilon)$, which implies $h_j^+(\hat{r})+h_j^-(\hat{r})>0$ for $t\in (0,\epsilon)$ by Lemma \ref{g-Delaunay lemma}.
Note that $\alpha(\hat{r})$ is continuous for $t\in[0,\epsilon]$, smooth for $t\in(0,\epsilon)$ and
		$$\frac{\partial \alpha}{\partial t}(\hat{r})=\frac{h_j^+(\hat{r})+h_j^-(\hat{r})}{l_{0j}}>0, t\in (0,\epsilon)$$
by Lemma \ref{basic property I of IDCP} (b),  we have $\alpha(\hat{r})>\alpha(r)$ for $t\in(0,\epsilon)$.
	\end{proof}
	Now we can prove Theorem \ref{Maximum principle}, which is paralleling to the proof of the maximal principle in \cite{LSW}.
For completeness, we include the proof here.

\textbf{Proof for Theorem \ref{Maximum principle}:}
Without loss of generality, we assume $\frac{r_0}{\bar{r_0}}=1$ and $r_i\le{\bar{r_i}}$ for all $i=1,2,\dots,n$,
otherwise we can scale
$r_i$ ($i\in \{0, \cdots, n\}$) by a factor $\frac{\overline{r}_0}{r_0}$. We prove the theorem by contradiction.
Otherwise, there exists a weighted Delaunay inversive distance circle packing $r$ on $(P_n, \mathcal{T}, I)$ such that $r_0=\bar{r}_0$, $r_i\le\bar{r}_i$ for all $i=1,2,\dots,n$ with one $r_{i_0}<\bar{r}_{i_0}$ and  $\alpha(\bar{r})\leq\alpha(r)$.
By Lemma \ref{key lemma},  after replacing $r$ by $\hat{r}$, we may assume that
\begin{equation}\label{alpha(bar-r)<alpha(r)}
  \alpha(\bar{r})<\alpha(r).
\end{equation}
		Consider the set
		\begin{equation*}
\begin{aligned}
		  X:=\{x\in\mathbb{R}^{n+1}|r\le x\le \bar{r}, &x \text{ is a generalized weighted Delaunay inversive} \\
&\text{distance circle packing on} (P_n, \mathcal{T}, I) \}.
\end{aligned}
		\end{equation*}
		Obviously, $r\in X$ and $X$ is bounded. By Lemma \ref{theta continuous}, $X$ is a closed set in $\mathbb{R}^{n+1}$.
Therefore, $X$ is a nonempty compact set and $\alpha(x)$ has a maximum point on $X$.
Let $t\in X$ be a maximum point of the continuous function $f(x)=\alpha(x)$ on $X$. If $t\ne\bar{r}$, then by Lemma \ref{key lemma},
we can find a weighted Delaunay inversive distance circle packing $\hat{t}$ on $(P_n, \mathcal{T}, I)$ such that $\hat{t}\ge t$, $\hat t_0=\bar{r}_0$, $\hat{t}\le \bar{r}$ and $\alpha(\hat{t})>\alpha(t)$, which implies that $t$ is not a maximum point of $f(x)=\alpha(x)$ on $X$. So $t=\bar{r}$ and then we have
		$$\alpha(\bar{r})=\alpha(t)\ge\alpha(r)>\alpha(\bar{r}),$$
where the last inequality comes from (\ref{alpha(bar-r)<alpha(r)}).
		This is a contradiction.
\qed

\begin{figure}[!htb]
\centering
\begin{overpic}[height=.2\textwidth,width=.4\textwidth]{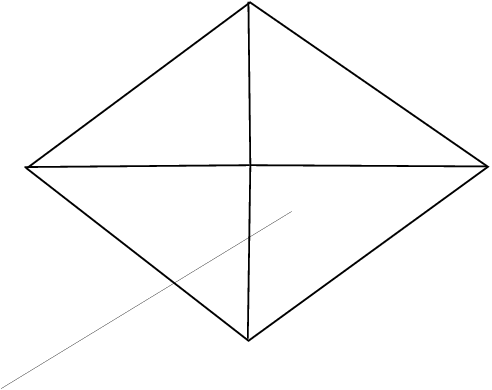}
  \put(50,0){$ v_{1}$}
    \put(52,32){$ v_{0}$}
  \put(50,52){$ v_{3}$}
  \put(0,31){$ v_{4}$}
  \put(100,31){$ v_{2}$}

\end{overpic}
  \caption{Counterexample for the maximal principle with $I$ in $(-1, +\infty)$}
  \label{counterexample}
\end{figure}

\begin{remark}
  The maximal principle is sharp in the sense that it can not be extended to the case that
  the weight $I$ takes value in $(-1, +\infty)$ and satisfies the structure condition (\ref{structure condition}). Specially, it does not allow the weight to take value in $(-1, 0)$ and $(1, +\infty)$ at the same time.
  We have the following counterexample for this case. Let $P_4$ be a polygon disk with four boundary vertices $v_1, v_2, v_3, v_4$ and a unique interior vertex $v_0$. Please refer to Figure \ref{counterexample}.
  Set $I_{01}=I_{03}=-\frac{1}{2}$, $I_{02}=I_{04}=2$ and $I_{12}=I_{23}=I_{34}=I_{41}=1$.
  Such a weight is regular and satisfies the structure condition (\ref{structure condition}).
  We further set $r_0=1$, $r_1=r_3=2$ and $r_2=r_4=c^{-1}$ with $c>0$.
  It is direct to check that, for each triangle $\triangle v_0v_iv_{i+1}, i=1,2,3,4$, $Q=\frac{3}{4}(c^2+4c+1)>0$ for any $c>0$.
  This implies the triangles $\triangle v_0v_iv_{i+1}$ are all non-degenerate and congruent. Furthermore, it is direct to check that
  $h_{1}^{+}(r)=h_{1}^{-}(r)=h_{3}^{+}(r)=h_{3}^{-}(r)>0$ and $h_{2}^{+}(r)=h_{2}^{-}(r)=h_{4}^{+}(r)=h_{4}^{-}(r)=0$ for any $c>0$. Therefore, $r=(1, 2, c^{-1}, 2, c^{-1})\in \mathbb{R}^5_{>0}$ is a non-degenerate weighted Delaunay inversive distance circle packing on $P_{4}$ for any $c>0$. We can also check that $l_{0i}^2+l_{0,i+1}^2=l_{i,i+1}^2$ for any $c>0$, which implies that the triangles $\triangle v_0v_iv_{i+1}$ are right triangles with $\angle v_iv_0v_{i+1}=\frac{\pi}{2}$.
  Therefore, the cone angle $\alpha(r)$ at the vertex $v_0$ is always $2\pi$ for any $c>0$. This implies that the maximal principle is not valid in this case.
\end{remark}

\subsection{A ring lemma}
\begin{lemma}\label{ring lemma}
Let $\mathcal{T}_{st}$ be the standard hexagonal triangulation of the plane and
$I:E\to (-\frac{1}{2}, +\infty)$ be a weigh defined on the edges.
$r:V \to (0, +\infty)$ is an inversive distance circle packing so that $(\mathcal{T}, I, r)$ is a geometric triangulation of the plane.  If $v_0\in V$, then for any $r:V\to (0, +\infty)$, there exists $C = C(v_0, I, \mathcal{T})>0$ such that
$$r(v_0) \leq Cr(v_k) \quad \text{if } v_k\in N(v_0).$$
\end{lemma}

\begin{figure}[!htb]
\centering
  \includegraphics[height=.3\textwidth,width=.45\textwidth]{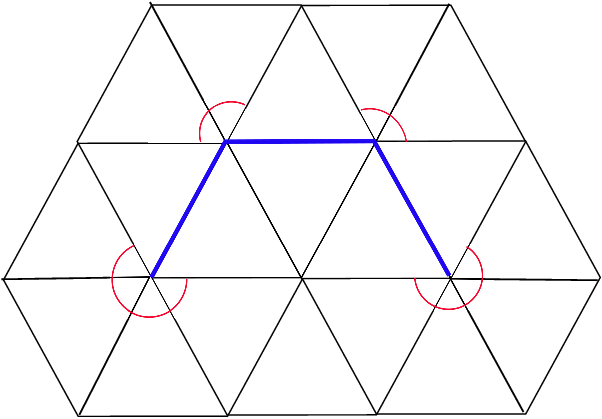}
  \caption{Ring Lemma.}
  \label{ringlemma}
\end{figure}

\begin{proof}
If not, we can assume that there exists a sequence of inversive distance circle packings $r_n: V\rightarrow (0, +\infty)$ such that $(T, I, r_n)$ is a geometric triangulation of the plane with $\lim_{n\to\infty} r_n(v_0) = 1$ and $\lim_{n\to\infty} r_n(v_1) = 0$ for $v_1\in N(v_0)$.
Here $N(v_0)$ denotes the vertices in $V$ adjacent to $v_0$.
By taking subsequences of $\{r_n\}$, we can assume that $r_n(v)$ converges in $[0, +\infty]$ for any $v\in V$. If $\lim_{n\to\infty} r_n(v) = 0$, then we call the vertex $v$ is degenerate.

Let $\mathcal{C}$ be the connected subcomplex of $\mathcal{T}$ generated by degenerate vertices such that $v_1\in \mathcal{C}$, and let $\mathcal{B}$ be the maximal connected subcomplex generated by vertices adjacent to vertices in $\mathcal{C}$.
Note that vertices in $\mathcal{B}$ are not isolated, otherwise the curvature at the vertex could not be zero.

We claim that there are at most five edges in $\mathcal{B}$ whose link intersects with $\mathcal{C}$. These edges are in the boundary of $\mathcal{B}$.
Otherwise, there are six triangles with one degenerate vertex and two non-degenerate vertices as $n\to \infty$. Note that by definition, the degenerate vertices are mapped to one point $O$ in the plane as $n\to \infty$. Hence, there are six triangles, each of which has one vertex mapped to $O$. By the assumption that $I>-1/2$, the angle of the degenerate vertices in these triangles are strictly larger than $\pi/3$. This implies that the curvature of $O$ can not be zero, and interiors of these six triangles are not disjoint from each other. This contradicts the fact that $(T, I, r_n)$ are geometric triangulations of the plane. This completes the proof of the claim.

Since the smallest cycle in $\mathcal{T}$ separating points has length six, these five edges (or fewer) can not form a loop which separates points. Then $\mathcal{B}$ is contractible, and all the vertices adjacent to vertices in $\mathcal{B}$ are degenerate by the maximality of $\mathcal{B}$. Then it is straightforward to check that the sum of the curvatures of vertices in $\mathcal{B}$ can not be zero. For example, if one connected component is $P_4$ shown in the Figure \ref{ringlemma}, then
$$\sum_{v\in P_4}K(v) = \sum_{v\in P_4} (2\pi - \sum_{(v,f)} \theta(v, f)),$$
where $(v, f)$ means that $f$ is a face in $\mathcal{T}$ containing $v$. Then since vertices adjacient to $\mathcal{B}$ are degenerate, then angles at $v$ are zero if the triangle $f$ containing $v$ has two degenerate vertices, as the red angles shown in Figure \ref{ringlemma}. Then there are at most six triangles containing the edges of $\mathcal{B}$, which contributes to the angle sum at vertices of $\mathcal{B}$. Therefore, we have
$$\sum_{v\in P_4}K(v) = \sum_{v\in P_4} (2\pi - \sum_{(v,f)} \theta(v, f)) \geq 8\pi - 6\times \pi = 2\pi >0.$$
This completes the proof.
\end{proof}

Lemma \ref{ring lemma} corresponds to the ring lemma in Rodin-Sullivan's famous work \cite{RS}.
Lemma \ref{ring lemma} depends on the hexagonal triangulation of the plane and the idea of its proof comes from \cite{H2}.
Following He's method in \cite{H2}, one can also prove a ring lemma for inversive distance circle packings with $I\in [0, +\infty)$ on surfaces with arbitrary triangulations.
Notice that the constant $C$ in Lemma \ref{ring lemma} depends on the vertex $v_0$ in general. However, if the weight comes from a lattice, then the constant $C$ works for all the vertices by periodicity.

\begin{corollary}\label{bounds on gradient}
Let $(\C, \T, w*l_{st})$ be a flat hexagonal triangulation of the plane discrete conformal to the hexagonal triangulation from a lattice with weight $I>-1/2$. Then there exists $M = M(I)>0$ such that
$$\sup_{i\sim j}|w_{i} - w_{j}| \leq M.$$
\end{corollary}

\subsection{Spiral hexagonal triangulations and linear discrete conformal factors}
We first recall the definition of developing maps in \cite{LSW}.
Let  $l$ be a flat polyhedral metric on a simply connected triangulated surface $(S, \mathcal{T})$. Its developing map $\phi:(S, \mathcal{T}, l)\to \mathbb{C}$ can be constructed by induction. We can start with any isometric embedding of a Euclidean triangle $t\in F$ in $\mathbb{C}$. This defines an initial map $\phi|t: (t, l)\to \mathbb{C}$, which can be extended to any adjacent triangle $s\in F$ such that $e = t\cap s \in E$ by isometrically embedding $s$ in $\mathbb{C}$ such that $\phi(e) = \phi(s)\cap \phi(t)$. Since $S$ is simply connected, we can continue this extension, which induces a well-defined developing map up to isometries of the plane.

\begin{proposition}\label{spiral}
Let $\mathcal{T}_{st}$ be the standard hexagonal triangulation of the plane and $I: E\rightarrow(-1, +\infty)$
be a weight defined on the edges satisfying the structure condition (\ref{structure condition}).
Let $l$ be a weighted Delaunay inversive distance circle packing metric determined by a label $u:V\rightarrow\mathbb{R}$ on  $(\mathbb{C}, \mathcal{T}_{st}, I)$ with the vertex set being a lattice $V= L = \{m\vec v_1 + n\vec v_2\}$, where $\{\vec v_1,\vec v_2\}$ is a geometric basis of the lattice $L$.
Suppose $w:V\to\mathbb{R}$ is a nonconstant linear function defined by two positive numbers $\lambda$ and $\mu$ via
\begin{equation}\label{expression of linear w}
  w(m\vec v_1 + n\vec v_2) =  m\log \lambda + n\log \mu
\end{equation}
and $w*l$ is a generalized weighted Delaunay inversive distance circle packing metric on $(\mathbb{C}, \mathcal{T}_{st}, I)$.
Then the following statements hold.
\begin{description}
  \item[(a)] $(\mathbb{C}, \mathcal{T}_{st}, I, w*l)$ is flat.
  \item[(b)] Let $\phi$ be the developing map for $(\mathbb{C}, \mathcal{T}_{st}, I, w*l)$. If there exists a non-degenerate triangle in $(\mathbb{C}, \mathcal{T}_{st}, I, w*l)$, then there are two different non-degenerate triangles $t_1$ and $t_2$ in $(\mathbb{C}, \mathcal{T}_{st}, I, w*l)$ such that $\phi(int(t_1))\cap \phi(int(t_2)) \neq \emptyset$. In other words, $\phi$ does not produce an embedding of $(\mathbb{C}, \mathcal{T}_{st}, I, w*l)$ in the plane.
  \item[(c)] If all the triangles in $(\mathbb{C}, \mathcal{T}_{st}, I, w*l)$ are degenerate, then there exists an automorphism $\psi$
  of the triangulation $\mathcal{T}_{st}$ and two constants $\gamma_1=\gamma_1(I, \vec v_1,\vec v_2)$ and $\gamma_2=\gamma_2(I, \vec v_1, \vec v_2)$ such that $w(\psi(m\vec v_1 + n\vec v_2))=m\ln \gamma_1 +n\ln \gamma_2$.
\end{description}
\end{proposition}
\begin{proof}
The proof of part (a) and (b) are the same as the proof for Proposition 3.4 in \cite{LSW}, so we omit the proof of part (a) and (b) here. We only present the proof of part (c).

To see part (c), since all the triangles in $(\mathbb{C}, \mathcal{T}_{st}, I, w*l)$ are degenerate, the inner angles of the triangles $t_1$ and $t_2$
are $0$ or $\pi$. Composing with an automorphism of the triangulation $\mathcal{T}_{st}$,
we may assume $\alpha_1=\gamma_2=\pi$, where
the angles are marked in Figure \ref{spiral hexagonal triangulation}.
	% \begin{figure}[h]
	% 	\includegraphics[width=0.4\textwidth]{spiral.png}
	% 	\caption{Angles of spiral hexagonal triangulations}
	% 		\label{spiral hexagonal triangulation}
	% \end{figure}
    \begin{figure}[h]
\begin{overpic}[scale=0.7]{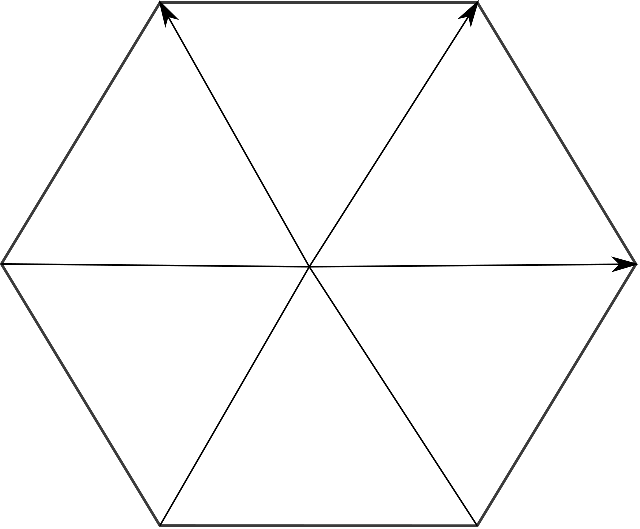}
 \put(5,80){$\vec v_{2} - \vec v_1$}

  \put(12,0){$-\vec v_{2}$}
  \put(80,0){$\vec v_{1}-\vec v_{2}$}

 \put(80,80){$\vec v_{2}$}
  \put(102,45){$\vec v_1$}
  \put(55,43){$\alpha_1$}
  \put(55,36){$\alpha_2$}
  \put(47,33){$\gamma_1$}
  \put(38,36){$\beta_2$}
    \put(90,36){$\beta_2$}
  \put(47,15){$t_1$}
  \put(47,70){$t_2$}
  \put(72,55){$t_1$}
  \put(72,25){$t_2$}
  \put(25,55){$t_1$}
  \put(25,25){$t_1$}

  \put(-10,45){$-\vec v_1$}

  \put(38,43){$\beta_1$}
  \put(47,48){$\gamma_2$}
  \put(90,43){$\beta_1$}
  \put(72,72){$\gamma_1$}
  \put(72,12){$\gamma_2$}

    \put(30,3){$\alpha_1$}
    \put(22,7){$\gamma_2$}
\end{overpic}
\caption{Angles of spiral triangulations.}
    \label{spiral hexagonal triangulation}
	\end{figure}
For the degenerate triangle with vertices $0$, $-\vec v_1$ and $-\vec v_2$, it is flat at $-\vec v_2$ by assumption, which implies $I_{0,-\vec v_1}>1$ by Lemma \ref{basic property I of IDCP} (c) and then $\Delta_{0,-\vec v_1,-\vec v_2}>0$ by Lemma \ref{basic property I of IDCP} (d). By Lemma \ref{simply connect of admi space with weight}, we further have
\begin{equation}\label{kappa* equation 1}
\begin{aligned}
  \kappa^*(-\vec v_2)=&\frac{1}{I_{0,-\vec v_1}^2-1}\{\gamma_{-\vec v_1,-\vec v_2,0}\kappa^*(0)+\gamma_{0, -\vec v_1,-\vec v_2}\kappa^*(-\vec v_1)\\
  &+\sqrt{\Delta_{0,-\vec v_1,-\vec v_2}[(\kappa^*(0))^2+(\kappa^*(-\vec v_1))^2+2I_{0,-\vec v_1}\kappa^*(0)\kappa^*(-\vec v_1)]}\},
  \end{aligned}
\end{equation}
where $\gamma_{v_i,v_j,v_k}:=I_{v_jv_k}+I_{v_iv_k}I_{v_iv_j}\geq 0$ by the structure condition (\ref{structure condition}) and
we use $^*$ to denote that we are discussing in the metric $w*l$.
Note that $\kappa^*(0)=\kappa(0)$, $\kappa^*(-\vec v_1)=\kappa(-\vec v_1)\lambda$ and
$\kappa^*(-\vec v_2)=\kappa(-\vec v_2)\mu$, we have
\begin{equation}\label{lambda-mu equation 1}
\begin{aligned}
\kappa(-\vec v_2)\mu=&\frac{1}{I_{0,-\vec v_1}^2-1}\{\gamma_{-\vec v_1,-\vec v_2,0}\kappa(0)+\gamma_{0, -\vec v_1,-\vec v_2}\kappa(-\vec v_1)\lambda\\
  &+\sqrt{\Delta_{0,-\vec v_1,-\vec v_2}[(\kappa(0))^2+\kappa^2(-\vec v_1)\lambda^2+2I_{0,-\vec v_1}\kappa(0)\kappa(-\vec v_1)\lambda]}\}
  \end{aligned}
\end{equation}
by (\ref{kappa* equation 1}).
Denote the right hand side of the equation (\ref{lambda-mu equation 1}) by $f_1(\lambda)$. Then $f_1(\lambda)$ is a strictly increasing function of $\lambda$
by $I_{0,-\vec v_1}>1$, $\Delta_{0,-\vec v_1,-\vec v_2}>0$ and the structure condition (\ref{structure condition}).
Furthermore, we have $\lim_{\lambda\rightarrow 0+}f_1(\lambda)=C_1>0$ and $\lim_{\lambda\rightarrow +\infty}f_1(\lambda)=+\infty$.
Dividing both sides of (\ref{lambda-mu equation 1}) by $\lambda$ gives
\begin{equation}\label{lambda-mu equation 2}
\begin{aligned}
\kappa(-\vec v_2)\frac{\mu}{\lambda}=&\frac{1}{I_{0,-\vec v_1}^2-1}\{\gamma_{-\vec v_1,-\vec v_2,0}\kappa(0)\lambda^{-1}+\gamma_{0, -\vec v_1,-\vec v_2}\kappa(-\vec v_1)\\
  &+\sqrt{\Delta_{0,-\vec v_1,-\vec v_2}[(\kappa(0))^2\lambda^{-2}+\kappa^2(-\vec v_1)+2I_{0,-\vec v_1}\kappa(0)\kappa(-\vec v_1)\lambda^{-1}]}\},
  \end{aligned}
\end{equation}
which implies that $\frac{\mu}{\lambda}$ is a strictly decreasing function of $\lambda$ with $\lim_{\lambda\rightarrow 0+}\frac{\mu}{\lambda}=+\infty$ and $\lim_{\lambda\rightarrow +\infty}\frac{\mu}{\lambda}=C_2>0$.

On the other hand, for the triangle with vertices $0$, $-\vec v_2$ and $\vec v_1-\vec v_2$, it is flat at $-\vec v_2$ by assumption,
which implies
$I_{0,\vec v_1-\vec v_2}>1$ by Lemma \ref{basic property I of IDCP} (c) and then $\Delta_{0,-\vec v_2,\vec v_1-\vec v_2}>0$ by Lemma \ref{basic property I of IDCP} (d). Applying Lemma \ref{simply connect of admi space with weight} to this triangle gives
\begin{equation}\label{kappa* equation 2}
\begin{aligned}
  \kappa^*(-\vec v_2)
  =&\frac{1}{I_{0,\vec v_1-\vec v_2}^2-1}\{\gamma_{0,-\vec v_2,\vec v_1-v_2}\kappa^*(\vec v_1-\vec v_2)+\gamma_{\vec v_1-\vec v_2, 0,-\vec v_2}\kappa^*(0)\\
  &+\sqrt{\Delta_{0,-\vec v_2,\vec v_1-\vec v_2}[(\kappa^*(0))^2+(\kappa^*(\vec v_1-\vec v_2))^2+2I_{0,-\vec v_1}\kappa^*(0)\kappa^*(\vec v_1-\vec v_2)]}\}.
  \end{aligned}
\end{equation}
Note that $\kappa^*(0)=\kappa(0)$,
$\kappa^*(-\vec v_2)=\kappa(-\vec v_2)\mu$ and $\kappa^*(\vec v_1-\vec v_2)=\kappa(\vec v_1-\vec v_2)\frac{\mu}{\lambda}$, we have
\begin{equation}\label{lambda-mu equation 3}
\begin{aligned}
  \kappa(-\vec v_2)\mu
  =&\frac{1}{I_{0,\vec v_1-\vec v_2}^2-1}\{\gamma_{0,-\vec v_2,\vec v_1-\vec v_2}\kappa(\vec v_1-\vec v_2)\frac{\mu}{\lambda}+\gamma_{\vec v_1-\vec v_2, 0,-\vec v_2}\kappa(0)\\
  &+\sqrt{\Delta_{0,-\vec v_2,\vec v_1-\vec v_2}[(\kappa(0))^2+\kappa^2(\vec v_1-\vec v_2)\frac{\mu^2}{\lambda^2}
  +2I_{0,-\vec v_1}\kappa^*(0)\kappa(\vec v_1-\vec v_2)\frac{\mu}{\lambda}]}\}
  \end{aligned}
\end{equation}
by (\ref{kappa* equation 2}).
Denote the right hand side of the equation (\ref{lambda-mu equation 3}) as $f_2(\lambda)$. Then $f_2(\lambda)$ is a strictly decreasing function
of $\lambda$ by $I_{0,\vec v_1-\vec v_2}>1$, $\Delta_{0,-\vec v_2,\vec v_1-\vec v_2}>0$, the structure condition (\ref{structure condition}) and
 the fact that $\frac{\mu}{\lambda}$ is a strictly decreasing function of $\lambda$.
Furthermore, $\lim_{\lambda\rightarrow 0+}f_2(\lambda)=+\infty$ and $\lim_{\lambda\rightarrow +\infty}f_2(\lambda)=C_3>0$.
Set $f(\lambda)=f_1(\lambda)-f_2(\lambda)$, then $f(\lambda)$ is a strictly increasing continuous function of $\lambda\in (0,+\infty)$ with
$\lim_{\lambda\rightarrow 0+}f(\lambda)=-\infty$ and $\lim_{\lambda\rightarrow +\infty}f_2(\lambda)=+\infty$,
which implies that there exists a unique number $\lambda=\lambda(I, \vec v_1, \vec v_2)\in (0,+\infty)$ such that $f_1(\lambda)=f_2(\lambda)$.
As a result, the system (\ref{lambda-mu equation 1}) and (\ref{lambda-mu equation 3}) has a unique solution $\lambda=\lambda(I, \vec v_1, \vec v_2)$ and $\mu=\mu(I, \vec v_1, \vec v_2)$ in $(0,+\infty)$. This completes the proof for part (c).
\end{proof}

\section{Rigidity of hexagonal triangulations of the plane}\label{section 4}
%\subsection{Rigidity of inversive distance circle packing.}
Recall the following definition and properties of embeddable flat polyhedral metrics in \cite{LSW}.

\begin{definition}[\cite{LSW} Definition 4.1]
Suppose $(S,\T)$ is a simply connected triangulated surface with a generalized PL metric $l$
and $\phi$ is developing map for $(S,\T, l)$.
$(S,\T, l,\phi)$ is said to be \it embeddable
\rm into $\C$ if for every simply connected finite subcomplex $P$
of $\T$, there exist a sequence of flat PL metrics on $P$ whose
developing maps $\phi_n: P \to \C$ are topological embeddings and converge uniformly to $\phi|_P$.
\end{definition}

\begin{lemma}[\cite{LSW} Lemma 4.2]\label{embed}
Let $(S,\T,l)$ be a flat polyhedral metric on
a simply connected surface with a developing map $\phi$.
\begin{enumerate}
    \item  Suppose $\phi$ is embeddable. If two simplices $s_1$, $s_2$ represent two distinct
non-degenerate triangles or two distinct edges in $\T$, then
$\phi(int(s_1))\cap \phi(int(s_2)) =\emptyset$.
    \item If $\phi$ is the pointwise convergent limit $\lim_{n\to
\infty}\psi_n$ of the developing maps $\psi_n$ of embeddable flat
 polyhedral metrics $(X, \T, l_n)$, then $(X,\T, l)$ is
embeddable.
\end{enumerate}
\end{lemma}

The standard hexagonal geodesic triangulations of open sets in $\C$ are embeddable. On the other hand, the generic Doyle spirals produce circle packings with overlapping interior, so the corresponding polyhedral metrics are not embeddable.

\begin{lemma}\label{flatness preserved by group action}
  Suppose $l_0$ is a weighted Delaunay inversive distance circle packing metric on $(S, \mathcal{T}_{st}, I)$ with
  the vertex set being a lattice $L=V$,
  the regular weight $I$ with values in $(-1, 1]$ or $[0,+\infty)$ satisfies
  the structure condition (\ref{structure condition}) and
   $I(e)=I(e+\delta)$ for any $e\in E$ and $\delta\in V$,
and  $l_0$ is generated by a label $w_0: V\rightarrow \mathbb{R}$.
  Suppose $(w-w_0)*l_0$ is a flat generalized weighted Delaunay inversive distance circle packing metric
  on the plane $(S, \mathcal{T}_{st})$.
  For any $\delta\in V$, set $u(v)=w(v+\delta)-w(v)$.
  Then $u*((w-w_0)*l_0)=(u+w-w_0)*l_0$ is a flat generalized weighed Delaunay inversive distance circle packing metric on $(S, \mathcal{T}_{st})$.
  Furthermore, if $u(v_0)=\max_{v\in V}u(v)$, then $u$ is a constant.
\end{lemma}

Lemma \ref{flatness preserved by group action} is a corollary of Theorem \ref{Maximum principle}, we omit the proof here.

\begin{lemma}\label{subsequence with const difference}
  Suppose $l_0$ is a weighted Delaunay inversive distance circle packing metric on $(S, \mathcal{T}_{st}, I)$ with
  the vertex set being a lattice $L=V$, the regular weight
  $I$ with values in $(-\frac{1}{2}, 1]$ or $[0,+\infty)$ satisfies  the structure condition (\ref{structure condition}) and
   $I(e)=I(e+\delta)$ for any $e\in E$ and $\delta\in V$,
and  $l_0$ is generated by a label $w_0: V\rightarrow \mathbb{R}$.
  Suppose $w*l_0$ is a flat generalized weighted Delaunay inversive distance circle packing metric
  on the plane $(S, \mathcal{T}_{st}, I)$.
  Then for any $\delta\in \{\pm u_1, \pm u_2, \pm(u_1-u_2)\}$, there exists $v_n\in V$ such that
  $$w_n(v):=w(v+v_n)+w_0(v+v_n)-w(v_n)-w_0(v_n)$$
  satisfies
\begin{description}
  \item[(a)] for all $v\in V$, the limit $w_\infty(v)=\lim_{n\rightarrow\infty}w_n(v)$ exists.
  \item[(b)] $(w_n-w_0)*l_0$ and $(w_\infty-w_0)*l_0$ are flat generalized weighted Delaunay inversive distance circle packing metrics
  on $(S, \mathcal{T}_{st}, I)$.
  \item[(c)] $w_\infty(v+\delta)-w_\infty(v)=a:=\sup \{w(v+\delta)-w(v)|v\in V\}$ for all $v\in V$.
  \item[(d)] the normalized developing maps $\phi_{(w_n-w_0)*l_0}$ of $(w_n-w_0)*l_0$ converges uniformly on compact subcomplex
  of $(S, \mathcal{T}_{st})$ to the normalized developing maps $\phi_{(w_\infty-w_0)*l_0}$ of $(w_\infty-w_0)*l_0$. As a result,
  if $(S, \mathcal{T}_{st}, I, w*l_0)$ is embeddable, then $(S, \mathcal{T}_{st}, I, (w_n-w_0)*l_0)$ is embeddable.
\end{description}
\end{lemma}
\proof
The proof is a modification of the proof of Lemma 4.5 in \cite{LSW}. For completeness, we include the proof here.
To see part (a), note that the label for the inversive distance circle packing metric $w*l_0$ is $w+w_0$.
By Lemma \ref{ring lemma}, there exists a constant
\begin{equation}\label{M}
\begin{aligned}
M=M(V, I)=\sup \{&|w(v+\delta)+w_0(v+\delta)-w(v)-w_0(v)||\\
&v\in V, \delta\in \{\pm u_1, \pm u_2, \pm(u_1-u_2)\}
\end{aligned}
\end{equation}
in $(0,+\infty)$ such that for fixed $\delta\in \{\pm u_1, \pm u_2, \pm(u_1-u_2)$,
$$a:=\sup \{w(v+\delta)+w_0(v+\delta)-w(v)-w_0(v)|v\in V\}\leq M.$$
Therefore, there exist a sequence $\{v_n\}$ in $V$ such that
\begin{equation}\label{w_n delta}
  a-\frac{1}{n}\leq w_n(\delta)=w(v_n+\delta)+w_0(v_n+\delta)-w(v_n)-w_0(v_n)\leq a.
\end{equation}
Furthermore, we have $w_n(0)=0$ and
\begin{equation}\label{w_n v+delta-w_n v}
w_n(v+\delta)-w_n(v)=w(v+\delta+v_n)+w_0(v+\delta+v_n)-w(v+v_n)-w_0(v+v_n)\leq a
\end{equation}
by the definition of $w_n$ and $a$.
By Lemma \ref{ring lemma}, if $v\in V$ is of combinatorial distance $m$ to $0$, then
\begin{equation*}
  \begin{aligned}
  |w_n(v)|=&|w_n(v)-w_n(0)|\\
  \leq& \sum_{i=1}^m |w_n(v_i)-w_n(v_{i-1})|\\
  =& \sum_{i=1}^m |w(v_i+v_n)+w_0(v_i+v_n)-w(v_{i-1}+v_n)-w_0(v_{i-1}+v_n)|\\
   \leq& m M
  \end{aligned}
\end{equation*}
by (\ref{M}),
where $v_m=v$, $v_0=0$ and $v_0\sim v_1\sim\cdots\sim v_m$ is a path of combinatorial distance $m$ between $0$ and $v$.
By the diagonal arguments, there exists a subsequence of $\{v_n\}$, still denoted by $\{v_n\}$ for simplicity, such that
$w_\infty(v):=\lim_{n\rightarrow \infty}w_n(v)$ exists for all $v\in V$.

To see part (b), for any fixed $n\in \mathbb{N}$ and any edge $e\in E$, we have
\begin{equation}\label{(w_n-w_0)*l_0}
  \begin{aligned}
(w_n-w_0)*l_0(e)=e^{-w(v_n)-w_0(v_n)}w*l_0(e+v_n)
  \end{aligned}
\end{equation}
by the translating invariance $I(e)=I(e+\delta)$ for the weight $I$. This implies that
$(w_n-w_0)*l_0$ is a flat generalized weighted Delaunay inversive distance circle packing metric on $(S, \mathcal{T}_{st}, I)$
by the assumption that $w*l_0$ is a flat generalized weighted Delaunay inversive distance circle packing metric
  on $(S, \mathcal{T}_{st}, I)$.
By $w_\infty(v)=\lim_{n\rightarrow \infty}w_n(v)$ and continuity, we have
$(w_\infty-w_0)*l_0$ is a flat generalized weighted Delaunay inversive distance circle packing metric on $(S, \mathcal{T}_{st}, I)$.
Similarly, we have $w_\infty(v+\delta)-w_\infty(v)\leq a$ for any $v\in V$ by (\ref{w_n v+delta-w_n v}),
which implies
\begin{equation}\label{w infty<a}
  \begin{aligned}\sup \{w_\infty(v+\delta)-w_\infty(v)|v\in V\}\leq a.
  \end{aligned}
\end{equation}

To see part (c), by $w_n(0)=0$, (\ref{w_n delta}) and (\ref{w infty<a}), we have $w_\infty(0)=0$ and
\begin{equation*}
  \begin{aligned}
w_\infty(\delta)-w_\infty(0)=w_\infty(\delta)=a\geq \sup \{w_\infty(v+\delta)-w_\infty(v)|v\in V\},
  \end{aligned}
\end{equation*}
which implies that $w_\infty(v+\delta)-w_\infty(v)$ attains the maximal value $\sup \{w_\infty(v+\delta)-w_\infty(v)|v\in V\}$ at $v=0$.
Note that for fixed $\delta$ and $u(v):=w_\infty(v+\delta)-w_\infty(v)$,
$u*((w_\infty-w_0)*l_0)$ is a flat generalized weighted Delaunay inversive distance circle packing metric on $(S, \mathcal{T}_{st}, I)$
by Lemma \ref{flatness preserved by group action}.
By the discrete maximal principle, i.e. Theorem \ref{Maximum principle}, we have $w_\infty(v+\delta)-w_\infty(v)=a$ for any $v\in V$.

If $(S, \mathcal{T}_{st}, I, w*l_0)$ is embeddable, then $(S, \mathcal{T}_{st}, I, (w_n-w_0)*l_0)$ is embeddable by (\ref{(w_n-w_0)*l_0}).
The rest of the proof is an application of Lemma \ref{embed}.
\qed

As a direct corollary of Lemma \ref{subsequence with const difference}, we have the following result in the case $w_0$ being a constant.
\begin{corollary}\label{subsequence with const difference for w0 const}
  Suppose $l_0$ is a weighted Delaunay inversive distance circle packing metric on $(S, \mathcal{T}_{st}, I)$ with
  the vertex set being a lattice $L=V$, the regular weight
  $I$ with values in $(-\frac{1}{2}, 1]$ or $[0,+\infty)$ satisfies the structure condition (\ref{structure condition}) and
   $I(e)=I(e+\delta)$ for any $e\in E$ and $\delta\in V$,
and  $l_0$ is generated by a constant label $w_0: V\rightarrow \mathbb{R}$.
  Suppose $w*l_0$ is a flat generalized weighted Delaunay inversive distance circle packing metric
  on the plane $(S, \mathcal{T}_{st}, I)$.
  Then for any $\delta\in \{\pm u_1, \pm u_2, \pm(u_1-u_2)\}$, there exists $v_n\in V$ such that
  $$w_n(v):=w(v+v_n)-w(v_n)$$
  satisfies
\begin{description}
  \item[(a)] for all $v\in V$, the limit $w_\infty(v)=\lim_{n\rightarrow\infty}w_n(v)$ exists.
  \item[(b)] $w_n*l_0$ and $w_\infty*l_0$ are flat generalized weighted Delaunay inversive distance circle packing metrics
  on $(S, \mathcal{T}_{st}, I)$.
  \item[(c)] $w_\infty(v+\delta)-w_\infty(v)=a:=\sup \{w(v+\delta)-w(v)|v\in V\}$ for all $v\in V$.
  \item[(d)] the normalized developing maps $\phi_{w_n*l_0}$ of $w_n*l_0$ converges uniformly on compact subcomplex
  of $(S, \mathcal{T}_{st})$ to the normalized developing maps $\phi_{w_\infty*l_0}$ of $w_\infty*l_0$. As a result,
  if $(S, \mathcal{T}_{st}, I, w*l_0)$ is embeddable, then $(S, \mathcal{T}_{st}, I, w_\infty*l_0)$ is embeddable.
\end{description}
\end{corollary}

\begin{remark}
As the weight $I$ satisfies the translating invariance $I(e)=I(e+\delta)$ in Lemma \ref{subsequence with const difference} and
Corollary \ref{subsequence with const difference for w0 const},
the weight $I$ is in fact determined by $I_{0u_1}, I_{0u_2}$, and  $I_{u_1u_2}$.
There are some further restrictions on the weight $I$ under the conditions in Corollary \ref{subsequence with const difference for w0 const}.
Consider the triangle $\triangle 0u_1u_2$, as $l_0$ is a weighted Delaunay inversive distance circle packing metric on $\triangle 0u_1u_2$,
we have
$$3-I^2_{1}-I^2_{2}-I^2_{3}+2I_1I_2+2I_1I_3+2I_2I_3+2I_1+2I_2+2I_3>0$$
by $w_0=const$ and Lemma \ref{basic property I of IDCP} (a), where $I_1=I_{0u_1}, I_2=I_{0u_2}, I_3=I_{u_1u_2}$.
Specially, $I_1=I_2=I_3\in [0,+\infty)$ satisfies the conditions on the weight $I$ in
Corollary \ref{subsequence with const difference for w0 const}.
\end{remark}

Theorem \ref{infrigidity introduction} is a special case of the following result.

\begin{theorem}
  Suppose $l_0$ is a weighted Delaunay inversive distance circle packing metric on $(S, \mathcal{T}_{st}, I)$ with
  the vertex set being a lattice $L=V$, the regular weight
  $I$ with values in $(-\frac{1}{2}, 1]$ or $[0,+\infty)$ satisfies   the structure condition (\ref{structure condition}) and
   $I(e)=I(e+\delta)$ for any $e\in E$ and $\delta\in V$,
and  $l_0$ is generated by a constant label $w_0: V\rightarrow \mathbb{R}$.
  Suppose $w*l_0$ is a flat generalized weighted Delaunay inversive distance circle packing metric
  on the plane $(S, \mathcal{T}_{st}, I)$ and $(S, \mathcal{T}_{st}, I,w*l_0)$ is embeddable into $\mathbb{C}$.
  Then $w$ is a constant function.
\end{theorem}

\begin{proof}
The idea of the proof follows the proof of Theorem 4.3 in \cite{LSW}.
We present the proof here for completeness. The idea can be summarized as follows. Assume $w$ is not a constant, we will construct a sequence of discrete conformal factor $w_n$ by extracting ``directional derivatives" of $w$ at different base points. This construction relies heavily on the symmetric structure of the lattice $V(\mathcal{T}_{st}) = L$ generated by $I$ and $w_0$, which implies that the limit of this sequence produce a linear discrete conformal factor $w_\infty$. By Corollary \ref{subsequence with const difference for w0 const}, $(S, \mathcal{T}_{st},I, w_\infty*l_0)$ is embeddable. However, by Proposition \ref{spiral}, if $w_\infty$ is not a constant,  $(S, \mathcal{T}_{st}, I, w_\infty*l_0)$ contains overlapping triangles under the developing maps. This leads to a contradiction.

%
%\textbf{Step 1}: directional derivative of $w$ satisfies the maximal principle.
%
%The first simple observation is that the symmetry of the lattice $L$ leads to a crucial properties for any discrete conformal factors of flat polyhedral metrics. Specifically, assume $\delta \in V$ and define $f: V \to \R$ by $f(v)=w(v+\delta) -w(v)$. It is immediate that $f*(w*l)=(f+w)*l = w(v+\delta)*l$ is a flat generalized
%weighted Delaunay inversive distance circle packing metric on $(S, \mathcal{T}_{st}, I)$.
%Then by the maximal principle, i.e. Theorem \ref{Maximum principle}, if there exists a vertex $v_0$ such that $f(v_0)=\max\{f(v) | v \in V\}$, then $f$ is constant.
%
%The function $f$ can be regarded as a discrete directional derivative of the discrete conformal factor $w$ in the direction of $\delta \in V$. The observation above means that the directional derivative $f$ of $w$ also satisfies the maximal principle.

\textbf{Step 1:} construct a linear limit $w_\infty$.
Since $w$ is assumed to be different from a constant function, then there exists a $\delta_1 \in L_0= \{\pm u_1, \pm u_2, \pm(u_1-u_2)\}$ such that $a_1 = \sup \{w(v + \delta_1) - w(v) | v\in V\}> 0$. By Corollary \ref{bounds on gradient}, $a_1 \in (0, \infty)$.
Applying Corollary \ref{subsequence with const difference for w0 const} to $w*l_0$ in the direction $\delta_1$,
there exists a function $w_\infty: \rightarrow \mathbb{R}$ such that $(S, \mathcal{T}_{st},I, w_\infty*l_0)$ is embeddable and
$$w_\infty(v+\delta_1)-w_\infty(v)=a_1, \forall v\in V.$$
Further applying Corollary \ref{subsequence with const difference for w0 const} to $w_\infty*l_0$ in the direction
$\delta_2\in \{\pm u_1, \pm u_2, \pm(u_1-u_2)\}-\{\pm \delta_1\}$ gives rise to a function
$\hat{w}=(w_\infty)_\infty: V\rightarrow \mathbb{R}$ such that $(S, \mathcal{T}_{st},I, \hat{w}*l_0)$ is embeddable and
$$\hat{w}(v+\delta_1)-\hat{w}(v+\delta_1)=a_1, \hat{w}(v+\delta_2)-\hat{w}(v)=a_2, \forall v\in V,$$
which shows that $\hat{w}(v)$ is an affine function of the form
$\hat{w}(n\delta_1+m\delta_2)=na_1+ma_2+a_3$ with $a_1\in (0, +\infty), a_2, a_3\in \mathbb{R}$.
Without loss of generality, we can assume $\hat{w}(n\delta_1+m\delta_2)=na_1+ma_2$ as the properties of weighted Delaunay and generalized PL metrics are invariant under scaling.
Then we obtain a function $\hat{w}: V\rightarrow \mathbb{R}$ satisfying
$\hat{w}(n\delta_1+m\delta_2)=na_1+ma_2$ and $(S, \mathcal{T}_{st},I, \hat{w}*l_0)$ is embeddable.

	\begin{figure}[h]
\begin{overpic}[angle=90, scale=0.6]{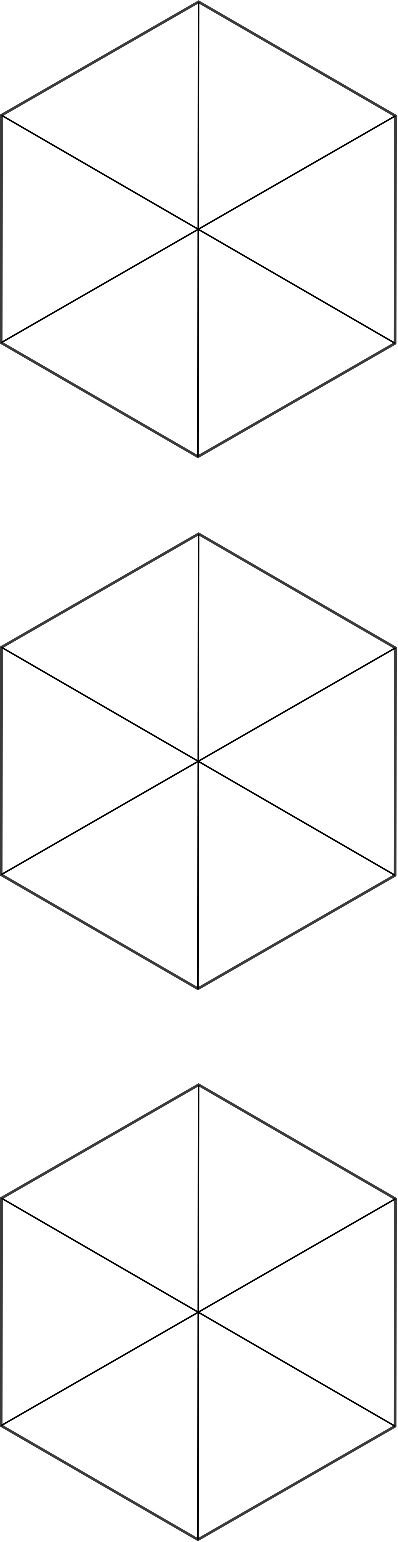}
 \put(21.5,22){$0$}
 \put(56,22){$0$}
 \put(92,22){$0$}

  \put(89,23){$\pi$}
 \put(53,23){$0$}
 \put(18.5,23){$0$}
 \put(10,23){$0$}
 \put(45,23){$\pi$}
 \put(80,23){$0$}

 \put(14,15){$\pi$}
  \put(49,15){$0$}
  \put(84.5,15){$0$}

\put(17,14){$\pi$}
\put(52,14){$\pi$}
\put(87,14){$\pi$}

\put(96,14){$0$}
\put(61,14){$0$}
\put(26,14){$0$}
\end{overpic}
\caption{Three cases of degenerate triangulations.}	
	\label{Three cases of degenerate triangles}
	
	\end{figure}

\textbf{Step 2:} Overlapping of $(S, \mathcal{T}_{st}, I, \hat{w}*l_0)$.

By step 1,  there are two positive numbers $\lambda\in (1,+\infty)$ and $\mu\in (0, +\infty)$ so that
$$\hat{w}(m\delta_1 + n\delta_2) = m\log \lambda + n\log \mu$$
and $(S, \mathcal{T}_{st},I, \hat{w}*l_0)$ is embeddable.
Then there is no non-degenerate triangle in the image of the developing map $\hat{\phi}$ for $(S, \mathcal{T}_{st}, I, \hat{w}*l_0)$, otherwise by Proposition \ref{spiral}, there are two triangles with overlapping interior.

Therefore, all the triangles in the image of $(S, \mathcal{T}_{st}, I, \hat{w})$ under $\hat{\phi}$ are degenerate.
All the angles are either $0$ and $\pi$.
There are three cases in Figure \ref{Three cases of degenerate triangles} showing triangles in the star of the origin. The last case can be ruled out by the weighted Delaunay condition. The first two cases are differed by a rotation $\gamma$. Therefore, we just need to consider Case 1.
By Proposition \ref{spiral}(c), the constants $\lambda$ and $\mu$ depend only on $I$ and $w_0$.

% 	\begin{figure}[h]
% 	\label{intersecting edges}
% 		\includegraphics[width=0.4\textwidth]{step3.png}
% 		\caption{Intersecting edges in the developing maps. }
% \label{Intersecting edges in the developing maps}
% 	\end{figure}

	\begin{figure}[h]
\begin{overpic}[scale=0.5]{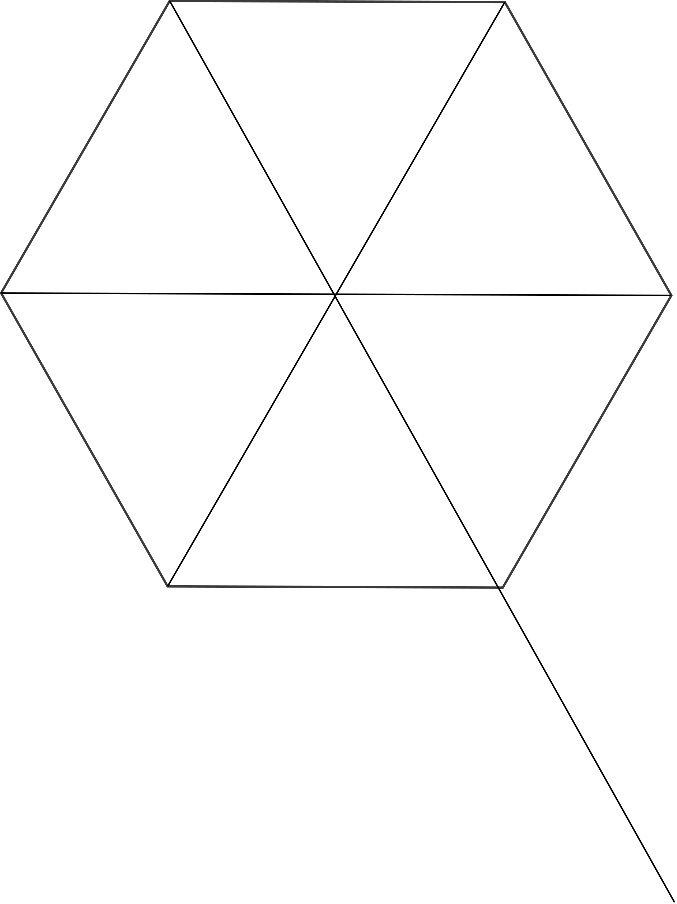}
 \put(42,70){$v_{0}$}
 \put(77,70){$v_{1}$}
  \put(55,102){$v_2$}
  \put(18,102){$v_3$}
   \put(5,70){$v_{4}$}
  \put(16,42){$v_5$}
  \put(53,42){$v_6$}
  \put(77,0){$v_7$}

\put(25, 80){$l_1$}
\put(48, 55){$l_2$}
\put(66, 20){$l_3$}
\end{overpic}
\caption{Intersecting edges in the developing maps. }
	\label{Intersecting edges in the developing maps}
	\end{figure}

Consider the lengths of edges $e_1 = v_0v_3$, $e_2 = v_0v_6$, and $e_3 = v_6v_7$ and their respective lengths $l_1$, $l_2$, and $l_3$ in $\hat{w}*l_0$ in Figure \ref{Intersecting edges in the developing maps}. Notice that $l_1 = (\lambda/\mu) l_2$ and $l_3 = (\mu/\lambda) l_2$, then
    $l_1 + l_3 \geq 2l_2 > l_2.$
Since  $(S, \mathcal{T}_{st}, I, \hat{w}*l_0)$ with the developing map $\hat{\phi}$ is embeddable, there exists a sequence of flat polyhedral metrics with developing maps $\phi_n$, which are embeddings, such that $\phi_n$ converges to $\hat{\phi}$ uniformly on compact sets. Then for $n$ large enough, the images of $e_1$ and $e_3$ under $\phi_n$ intersects by the inequality above.
This contradicts that $(S, \mathcal{T}_{st},I, \hat{w}*l_0)$ is embeddable, which completes the proof.
\end{proof}

\end{document}